\documentclass{article}
\usepackage{graphicx} 
\usepackage{latexsym} 
\usepackage{amsfonts} 
\usepackage{psfrag}
\usepackage{color} 

\newtheorem{theorem}{Theorem}[section]
\newtheorem{proposition}[theorem]{Proposition}
\newtheorem{definition}[theorem]{Definition}
\newtheorem{lemma}[theorem]{Lemma}
\newtheorem{corollary}[theorem]{Corollary}

\newcommand{\proof}{\noindent {\it Proof: }}
\newcommand{\qed}{\vrule height4pt width4pt depth4pt}

\newcommand{\beyond}[1]{{\mbox{\scriptsize beyond $#1$}}}
 
\def\endofproof{\hfill \nobreak \hskip0pt plus 1fill \qquad \qed
  \medskip\noindent}
 
\def\qed{\vbox{\hrule \hbox{\vrule\hbox to 5pt{\vbox to
        6pt{\vfil}\hfil}\vrule}\hrule}}

\newcommand{\nota}[1]{
}

\begin{document}

\pagestyle{headings}
\markboth{\hfill \today}{\today \hfill}
\addtolength{\oddsidemargin}{-2cm}
\addtolength{\evensidemargin}{-2cm}

\bibliographystyle{alpha}
\unitlength=1cm

\title{\sc The Complexity of Finding Small Triangulations of Convex $3$-Polytopes.}
\author{%
 Alexander Below%
\thanks{Institut f\"ur Theoretische Informatik, ETH-Z\"urich
({\tt below@inf.ethz.ch}).}
\and
Jes\'us A. De Loera%
\thanks{Dept. of Mathematics, Univ. of California-Davis
({\tt deloera@math.ucdavis.edu}).}
\and
J\"urgen Richter-Gebert%
\thanks{Institut f\"ur Theoretische Informatik, ETH-Z\"urich
({\tt richter@inf.ethz.ch}).}}

\maketitle

\begin{abstract}
  The problem of finding a triangulation of a convex three-dimensional
  polytope with few tetrahedra is $NP$-hard.  We discuss other related
  complexity results.
\end{abstract}

\section{Introduction}

A {\em triangulation} of a $d$-dimensional convex polytope $P$ is a
set of $d$-simplices whose union is the polytope, their vertices are
extreme points of $P$, and any two simplices in it intersect in a
common (possibly empty) face. The {\em size} of a triangulation is the
number of its full-dimensional simplices. In this paper we discuss the
computational complexity of finding small size triangulations of a
convex polytope. We discuss in particular the case of {\em minimal
triangulations}, i.e. those with smallest possible size.

This geometric minimization problem arises, for example, in Algebra
and Mathematical Programming. For example, minimal triangulations of
the $d$-cube have been extensively studied (see references in
\cite{cottle,haiman,hughes}) due to connections with the simplicial
approximation of fixed points of continuous maps (see \cite{todd}).
Optimal size triangulations appear also in the polyhedral techniques
in Algebraic Geometry \cite{surveybernd}. Understanding minimal
triangulations of convex polytopes is related to the problem of
characterizing the $f$-vectors of triangulations of balls and
polytopes (see open problems in \cite{handbookbilbjo}). In fact, the
study of minimal triangulations of topological balls also received
attention due to its connections to data structures, in the
calculation of rotation distance of binary trees
\cite{sleatoretal}. 

The computational geometry literature has several
papers interested in finding triangulations of optimal size
\cite{AvisE,edelsbrunner et al}. In 1992 Bern and Eppstein asked
whether there is a polynomial time algorithm to compute a minimal
triangulation of a 3-dimensional convex polytope (open problem 12 in
section 3.2 \cite{bern et al}). Our main result shows that, under the
hypothesis $P \not= NP$, such an algorithm cannot exist:

\begin{theorem} \label{main}
  Given a convex $3$-polytope $P$ with $n$ vertices and a positive
  integer $K$, deciding whether $P$ has a triangulation of size $K$ or
  less is an NP-complete problem.
\end{theorem}

We have the following corollaries (the second result was recently
obtained in \cite{Rich99} via a direct transformation to $3$-SAT):

\begin{corollary} \label{maincoro} 
  \begin{enumerate}
  \item Finding a minimal-size triangulation of a convex polytope of any
    fixed dimension $d\geq 3$ is NP-hard. Clearly, the same holds when the
    dimension is not fixed.

  \item Finding a minimal-size triangulation of the boundary of a convex
    polytope of any fixed dimension $d\geq 4$ is NP-hard. Clearly, the
    same holds when the dimension is not fixed.
  \end{enumerate}
\end{corollary}

\noindent Now we discuss the general structure and main ideas of 
the proof of Theorem~\ref{main}.  

We give a transformation to the Satisfiability (SAT) problem
(cf.~\cite{GarJohn79}): given an instance $S$ of $C$ logical clauses
in $V$ boolean variables, is there a truth assignment to the variables
such that all clauses are simultaneously satisfied?  We will give a
number $K$ and construct a convex 3-polytope, of size polynomial in
$C$ and $V$ (polynomial size pertains to the binary encoding length),
which has a triangulation of size at most $K$ if and only if there is
a satisfying truth assignment.  In fact we can restrict our discussion
to the special case of the SAT problem where each variable appears in
three clauses, two of the times negated (see page 259
\cite{GarJohn79}).

\begin{figure}[h]
  \psfrag{q0}{\small $q_0$}
  \psfrag{q1}{\small $q_1$}
  \psfrag{q2}{\small $q_2$}
  \psfrag{qm}{\small $q_m$}
  \psfrag{qmp1}{\small $q_{m+1}$}
  \psfrag{a}{\small $a$}
  \psfrag{b}{\small $b$}
  \centerline{
    \includegraphics[scale=.5]{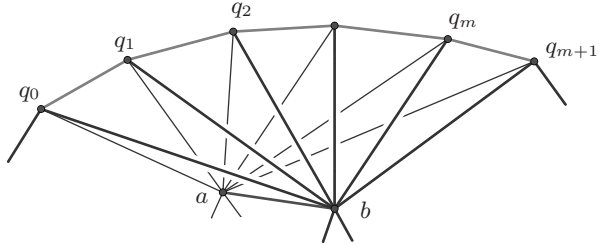} 
    }
  \caption{The vertex-edge chain of Lemma \ref{keyStructure}}
  \label{keyStructureFigure}
\end{figure}

Two elementary properties of triangulations will be useful to reach
our goal: (1) Every boundary triangular facet $F$ of a polytope is
contained in exactly one tetrahedron of a triangulation.  The fourth
vertex of that tetrahedron is said to {\em triangulate } $F$.  (2)
Simplices of a triangulation cannot intersect in their relative
interiors.  We will primarily see this behavior in triangles being
pierced by an edge of the triangulation, a so-called {\em bad
intersection}. Our proof combines techniques presented in the articles
\cite{Below99} and \cite{RupSei92}:

Consider a long vertex-edge chain on the boundary of a polytope such
that the adjacent triangular faces all meet in two points $a$ and $b$
(see Figure \ref{keyStructureFigure}).  The proof of the following
lemma can be found in Section 2 of \cite{Below99}:

\begin{lemma} \label{keyStructure}
  Let $P$ be a convex 3-polytope such that the triangles
  $(a,q_i,q_{i+1})$ and $(b,q_i,q_{i+1})$ for $i=0, \dots,m$ are among
  its facets, with the additional restriction that $\mbox{conv}\{a,b\}
  \cap \mbox{conv}\{q_0, \ldots, q_{m+1}\} = \emptyset$.  Let $n$ be
  the number of vertices of $P$.
  
  Then, for each triangulation of $P$ that does not use the (interior)
  edge $(a, b)$ the number of tetrahedra is at least $n + m - 3$.
\end{lemma}

For us, the number $n + m - 3$ will be a relatively large number, such
that a triangulation not using $(a, b)$ cannot be small. When using
the edge $(a, b)$ we can triangulate $P_{\mbox{\scriptsize chain}} =
\mbox{\textit{conv}}\{a, b, q_0, \ldots, q_{m+1}\}$ using the $m+1$
tetrahedra $(a, b, q_i, q_{i+1})$ for $i = 0, \ldots, m$.  Call $Q$
the (non-convex) polytope we get after cutting all these tetrahedra
out of $P$.  Let $n_Q$ denote the number of vertices of $Q$.  Note
that $n = n_Q + m$.  Suppose the number of tetrahedra in any
triangulation of $Q$ (if there is one at all) is bounded above by some
number $t$.  Then we can bound the size of a minimal triangulation of
$P$ {\em using} $(a,b)$ by $t + m$.  Note that $t$ does not depend on
the length $m$ of the vertex-edge chain.  Hence, by choosing $m$ large
enough (leaving $Q$ as it is) makes
\[
t + m < n_Q + 2m - 3,
\]
and any close-to-minimal triangulation will {\em have to} use $(a,
b)$.  This argument still holds when we have many vertex-edge chains
of the same length $m$ present in other parts of the boundary of the
polytope $P$.  If $m$ is large enough, a small triangulation
is forced to use the edges $(a, b)$ of {\em all} these vertex-edge chains.

\begin{figure}[h]
  \centerline{
    \psfrag{A1}{$A_1$}
    \psfrag{A2}{$A_2$}
    \psfrag{A3}{$A_3$}
    \psfrag{B1}{$B_1$}
    \psfrag{B2}{$B_2$}
    \psfrag{B3}{$B_3$}
    \includegraphics[scale=.6]{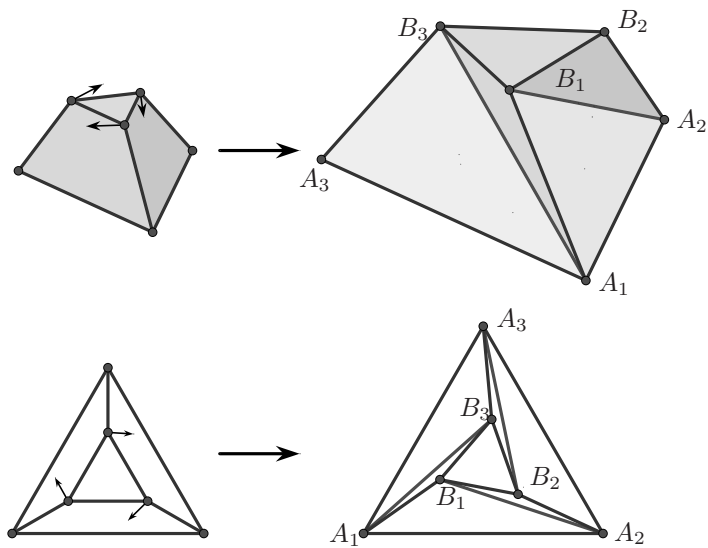}
    }
  \caption{A Sch\"onhardt polytope}
  \label{Schonhardt}
\end{figure}

We also use the famous non-convex {\em Sch\"onhardt polytope}
\cite{bern et al, handbooklee, lee, ORou87, RupSei92, Sch28}.  Roughly
speaking, a Sch\"onhardt polytope can be obtained by ``twisting'' the
top face of a triangular prism in a clockwise direction (see Figure
\ref{Schonhardt}).  The three quadrangular sides are then broken up
and ``bent in'', thus creating the non-convex ({\em reflex}) edges
$(B_i, A_{i+1})$ that we call {\em diagonals}.  The resulting polytope
is non-convex and we distinguish two triangular faces {\em the bottom}
$(A_1,A_2,A_3)$, and the {\em top or skylight} $(B_1,B_2,B_3)$ not
having a reflex edge (Note: whenever dealing with vertices of a
Sch\"onhardt polytope, abusing the notation, by an index $i + 1$ we
mean $(i \quad \mbox{mod} \quad 3) + 1$.  For example $3 + 1$ gives
$1$).

The reader can easily verify that Sch\"onhardt polytopes cannot be
triangulated using only its six vertices.  Imagine the Sch\"onhardt
polytope is glued along its bottom face to a bigger polytope.  Again,
the resulting non-convex polytope can only be triangulated if its top
triangle (the skylight) is {\em visible} by another vertex (we will
rigorously define visible later, but it does correspond to the
intuition that every point of the skylight can be joined to the vertex
by a straight line segment.  We will show that the triangular cone
defined by the hyperplanes containing the faces $(B_i, B_{i + 1}, A_{i
  + 1})$ contains exactly the points that can view the skylight.  For
this reason we will call this cone the {\em visibility cone}.
\begin{figure}
  \psfrag{visibilityCone}[r][r]{\small visibility cone}
  \psfrag{pillar}{\small diagonal}
  \psfrag{skylight}[c][c]{\small skylight}
  \psfrag{A1}{\small $A_1$}
  \psfrag{A2}{\small $A_2$}
  \psfrag{A3}{\small $A_3$}
  \psfrag{B1}{\small $B_1$}
  \psfrag{B2}{\small $B_2$}
  \psfrag{B3}{\small $B_3$}
  \hfill
  \includegraphics[scale=.5]{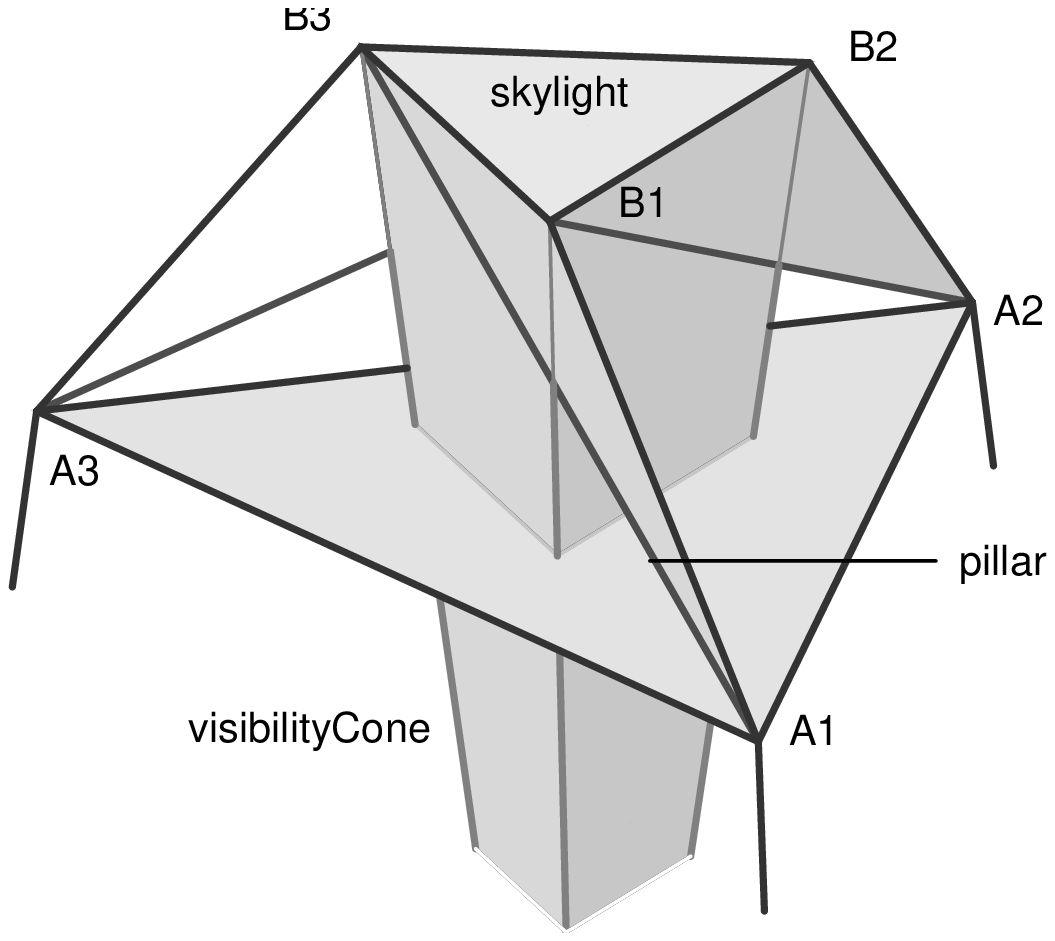}
  \hfill
  \includegraphics[scale=.4]{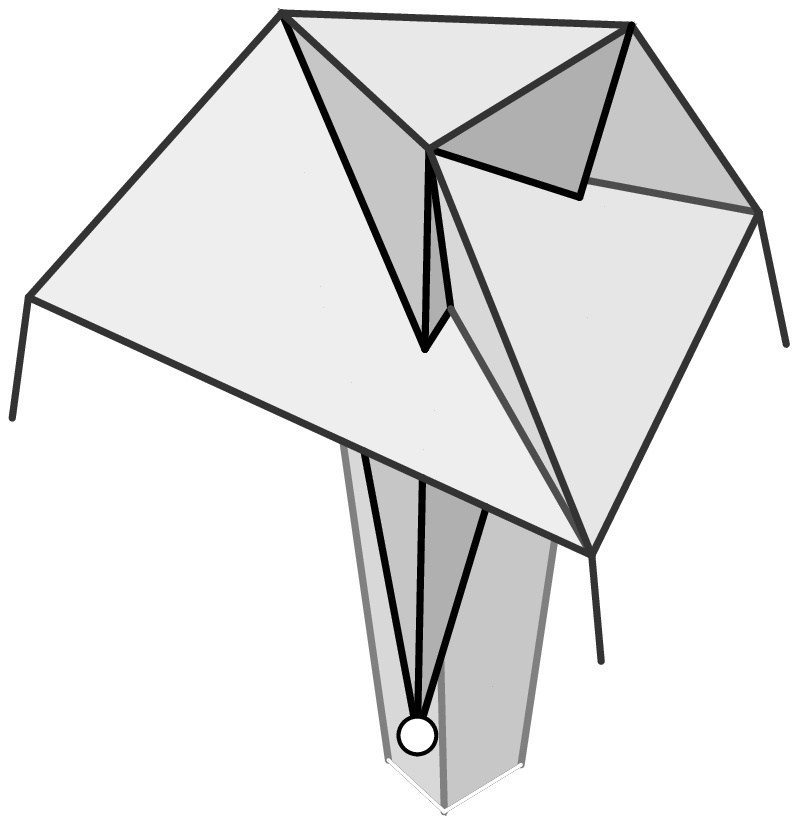}
  \hfill
  \caption{The visibility cone and an invisible vertex}
  \label{visibilityCone}
  \label{outsideVisibility}
\end{figure}

Now we convexify the Sch\"onhardt polytope by attaching three circular
vertex-edge chains opposite to the concavities.  This will give us a
convex polytope that we call the {\em cupola}, see Figure
\ref{introCupola}.
\begin{figure}[ht]
  \centerline{ \includegraphics[scale=.5]{13CupolaIntro.epsFinal} }
  \caption{Cupola}
  \label{introCupola}
\end{figure}
The cupola is usually glued along its bottom face to a bigger convex
polytope and obtain a convex polytope $P$.  We can combine what we
know about vertex-edge chains and about Sch\"onhardt polytopes.
Namely, in order to have a small triangulation of $P$, the three
diagonals of the Sch\"onhard polytope inside the cupola have to be
used.  But then, the vertex triangulating the skylight of the cupola
must not be obstructed from seeing the skylight by the diagonals.
Hence the vertex has to lie in the visibility cone of the cupola.

In \cite{RupSei92} Ruppert and Seidel used SAT to prove that it is
NP-complete to decide whether a non-convex polyhedron admits a
triangulation.  Their constructions used Sch\"onhardt polytopes, and
in particular their visibility cones, to do the transformation. In our
case, because we need convexity, we glue cupolas, instead of
Sch\"onhardt polytopes. They are glued to a bigger {\em frame}
polytope along their bottom faces.  Similar to \cite{RupSei92}, we
have variable cupolas and clause cupolas.  The visibility cones of the
variable cupolas contain only two {\em truth-setting} vertices, one
for {\em false} and one for {\em true}.  The visibility cones of the
clause cupolas contain as many {\em literal} vertices as there are
literals in the logical clause.  Each variable must choose between a
``true'' or ``false'' value. Inside each clause at least one variable
will be chosen to be true (to satisfy the clause). We model these
logical choices by the geometric choices of which vertex in the
visibility cone of a (variable/clause) cupola is used to triangulate
the skylight.  In addition, our polytope satisfies some {\em blocking}
conditions: the tetrahedron spanned by the top face of clause cupola
and a literal vertex coming from a negated variable $X_i$ will
improperly intersect the tetrahedron spanned by the top face of the
cupola of variable $X_i$ and the truth-setting vertex corresponding to
{\em true}.  In this way the choices made for the truth values of the
variables and for the literals satisfying the clauses will be {\em
  consistent}.  We will call our polytope the {\em logical polytope}
because it comes from a logical formula.

For the polynomial transformation (from SAT) we need to give an
algorithm to compute the coordinates of the logical polytope.  The
binary encoding length of the polytope, as well as the runtime of the
algorithm, have to be polynomial in the encoding length of the SAT
instance.  Each step of the construction will be polynomial, this is a
delicate point in the formalism of our argument.  We apply a sequence
of these constructions (polynomially many).  The coordinates of the
vertices of the polytope are potentially singly-exponential, but their
binary encoding length is guaranteed to be polynomial.

Elementary steps of construction include operations such as taking the
join of two or three points, intersecting planes and lines, putting
points on polynomial curves, etc.  The coordinates of the resulting
construction elements are therefore polynomials in coordinates of the
input elements.  On the other hand, we will have requirements on the
positions of the points with respect to some planes or other points on
lines etc.  All these conditions can be formulated as strict
polynomial inequalities in coordinates of the construction elements.
An essential element of our construction is that our systems of strict
polynomial inequalities will depend on {\em one single parameter}
$\epsilon$.  All these polynomial inequalities are satisfied at
$\epsilon = 0$, but an additional requirement for us is $\epsilon >
0$.  The following lemma describes a polynomial algorithm to find a
number $\epsilon_0$ such that all $0 < \epsilon \le \epsilon_0$ solve
the inequality system.
  
\begin{lemma} \label{openCondition} 
  \begin{enumerate} 
  \item Suppose $p(\epsilon) = a_d \epsilon^d + \cdots + a_1 \epsilon
    + a_0$ is a polynomial with $p(0) > 0$.  Let $\displaystyle
    \epsilon_0(p) := \min\left(1, \frac{\displaystyle
        a_0}{\displaystyle 2(|a_1| + \cdots + |a_d|)}\right)$.  Then
    for $0 \le \epsilon \le \epsilon_0(p)$ we have $p(\epsilon) > 0$.
    
    Hence, the construction of $\epsilon_0$ can be done in time
    polynomial in the encoding length of the coefficients of $p$, and
    $\epsilon_0$ has polynomial encoding length.
    
  \item $p_1, \ldots, p_l$ are univariate polynomials such that
    $p_1(0) > 0, \ldots, p_l(0) > 0$ then there is a rational number,
    $\epsilon_0 > 0$, such that $p_1(\epsilon) > 0, \ldots,
    p_l(\epsilon) > 0$ for all $0 < \epsilon \le \epsilon_0$.
    Moreover, the encoding length of $\epsilon_0$ is bounded by a
    polynomial in the encoding length of the coefficients of $p_1,
    \ldots, p_l$.
  \end{enumerate}
\end{lemma}
\proof For $0 \le \epsilon \le 1$ we have that $a_i \epsilon^i \ge
-|a_i| \epsilon$.  The reason is that for $a_i \ge 0$, $a_i \epsilon^i
\ge 0 \ge -|a_i| \epsilon$, and for $a_i < 0$, $a_i \epsilon^i > a_i
\epsilon = -|a_i| \epsilon$.  Hence for $0 \le \epsilon \le
\epsilon_0(p)$
\[
p(\epsilon) \ge \sum_{i = 1}^d -|a_i| \epsilon + a_0 > -\sum_{i =
  1}^d |a_i| \frac{a_0}{2 \sum_{i = 1}^d |a_i|} + a_0 > 0.
\]

For the second part, take the value $\epsilon_0(p_1, \ldots, p_r) :=
\min(\epsilon_0(p_1), \ldots, \epsilon_0(p_r)$).  Now {\em all} the
conditions are simultaneously satisfied. \endofproof
  
Of course, in general the real solutions of a multivariate system of
inequalities coming imposed by geometric requirements may be empty,
but our steps of construction reduce everything to sequentially
solving easy univariate systems of inequalities.

Here is the organization of our paper.  In Section \ref{illumination}
we discuss useful properties of Sch\"onhardt polytopes and of cupolas.
Later in the section we explain how to construct and glue cupolas that
have a prescribed visibility cone and how to construct visibility
cones that fit our purposes.  The polynomial transformation to SAT is
presented in Section \ref{transformation}.  From a given SAT instance
we construct a frame polytope to which we then glue the cupolas.  In
the final Section \ref{conclusions}, we discuss the consequences of
our result and related complexity problems.  As a complement of our
main theorem we present a family of polytopes (the so-called stacked
polytopes) for which the decision of Theorem \ref{main} can be solved
in polynomial time.

\section{Basic Building Blocks} \label{illumination}

\psfrag{A1}{\small $A_1$}
\psfrag{A2}{\small $A_2$}
\psfrag{A3}{\small $A_3$}
\psfrag{B1}{\small $B_1$}
\psfrag{B2}{\small $B_2$}
\psfrag{B3}{\small $B_3$}
\psfrag{P}{\small $P$}
\psfrag{H}{\small $H$}
\psfrag{r1}{\small $l_1$}
\psfrag{r2}{\small $l_2$}
\psfrag{r3}{\small $l_3$}
\psfrag{l1}{\small $l_1$}
\psfrag{l2}{\small $l_2$}
\psfrag{l3}{\small $l_3$}
\psfrag{D1}{\small $D_1$}
\psfrag{D2}{\small $D_2$}
\psfrag{D3}{\small $D_3$}
\psfrag{V}{\small $V$}
\psfrag{F}{\small $F$}
\psfrag{Rstar}[r][r]{\small $R_*$}
\psfrag{apex}[br][br]{\small 
  \begin{tabular}{r} apex of $V$, \\[-.1cm] also $R_*$ \end{tabular}
}
\psfrag{A1s}{\small $A_1'$}
\psfrag{A2s}{\small $A_2'$}
\psfrag{A3s}{\small $A_3'$}
\psfrag{PcapH}{\small $H \cap (P \setminus F)$}
\psfrag{Hs}{\small $H'$}

\sloppy

We recall the notion of {\em beyond a face} (see \cite{Zie94}): A
point $p$ is {\em beyond} a face $F$ of a polytope $P$ if it
(strictly) violates all inequalities defining facets of $P$ containing
$F$, but it strictly satisfies all other inequalities that define other
facets of $P$.  The polytope $P_\beyond{F}$ is the (closure of the )
set of all points beyond $F$.  We denote by $P \setminus F$ the
polyhedron defined by all facet--defining inequalities that do not
hold with equality in $F$.  This is exactly $P \cup P_\beyond{F}$.  In
our constructions we will often put one or more points beyond some
face, and then take the convex hull.  This will only destroy the
facets containing this face, and introduce new ones containing the new
points.  We will say we {\em attach} one polytope $P$ to another $Q$
{\em along} a facets $F_P$ of $P$ and $F_Q$ of $Q$ if $P \subseteq
Q_\beyond{F_Q}$ and $Q \subseteq P_\beyond{F_P}$.  It is important to
observe that the convex hull of their union contains both the face
lattices of $P$ and $Q$ without, of course, $F_P$ and $F_Q$.

\subsection{The Sch\"onhardt Polytope}

Let us turn to a well-known example of a non-convex non-triangulable
polytope, the so-called Sch\"onhardt polytope (named after its first
occurrence in \cite{Sch28}. See also \cite{ORou87}).  For the notion
of non-convex polytope and what it means to triangulate them we refer
to \cite{Chaz90}.

\begin{definition} A {\em Sch\"onhardt} polytope (Figure \ref{Schonhardt}) 
  is a non-convex polytope with six vertices $A_1$, $A_2$, $A_3$,
  $B_1$, $B_2$, and $B_3$ and facets $(A_1, A_2, A_3)$, $(B_1, B_2,
  B_3)$, $(A_1, B_1, A_2)$, $(B_1, A_2, B_2)$, $(A_2, B_2, A_3)$,
  $(B_2, A_3, B_3)$, $(A_3, B_3, B_1)$, and $(B_3, B_1, A_1)$.  At
  exactly the edges $(B_1, A_2)$, $(B_2, A_3)$, $(B_3, A_1)$ the
  corresponding facets are to span an interior angle greater than
  $\pi$ (the edges are said to be {\em reflex}).  These edges are
  called the {\em diagonals} of the Sch\"onhardt polytope.  The top face
  $(B_1, B_2, B_3)$ is called the {\em skylight} of the Sch\"onhardt
  polytope.
  
  Six points are said to be in {\em Sch\"onhardt position} if they are
  the vertices of a Sch\"onhardt polytope.  We say that the skylight
  is {\em visible} from a point $x$ (or $x$ is {\em able to see} the
  skylight, or $x$ is a {\em viewpoint} of the skylight) if the
  tetrahedron spanned by $x$ and the skylight does not intersect any
  of the diagonals in their relative interior. The {\em visibility
    cone} of the Sch\"onhardt polytope is the triangular cone bounded
  by the planes $B_1 B_2 A_2$, $B_2 B_3 A_3$, and $B_3 B_1 A_1$.  See
  Figure \ref{visibilityCone}.
\end{definition}

The use of the word ``skylight'' is motivated by the idea that the
skylight triangle is a glass window and light comes through it
illuminating the interior of the Sch\"onhardt polytope defining a cone
of light.  It is obvious that this non-convex polytope cannot be
triangulated (without adding new points): The fourth point of the
tetrahedron containing the skylight must be one of $A_1$, $A_2$, or
$A_3$, but the diagonals ``obstruct the view'' of the skylight from
these vertices.

It is our intention to patch the sides of the Sch\"onhardt polytope
with vertex-edge chains in order to convexify it (and then glue it to
a frame polytope). According to Lemma \ref{keyStructure}, a small
triangulation of this convex polytope must necessarily contain the
diagonals. In this case, the fourth point of the tetrahedron
containing the skylight also has to be able to see the skylight.  We
will show where to place the vertex-edge chains in order for them not
to be visible from the skylight.  Hence, the triangulating vertex has
to lie beyond $(A_1, A_2, A_3)$, we will show that it has to lie in
the visibility cone.

\begin{lemma} \label{SchonhardtPolytopeLemma}
  Let $A_1$, $A_2$, $A_3$, $B_1$, $B_2$, $B_3$ be six points in
  Sch\"onhardt position.  We denote by $C_{A,B}$ the convex hull of the
  six points. Then
  \begin{enumerate}
  \item All orientations of simplices spanned by four of these six
    points are determined up to one global sign change.  As a
    consequence, the six points are in convex position, and their
    convex hull $C_{A,B}$ is an octahedron
    that has $(A_1, A_2, A_3)$ and $(B_1, B_2, B_3)$ as facets
    and it has edges $(A_i, B_{i + 1})$ ($i = 1, 2, 3$).
  
\item There are no points that can see the skylight $(B_1, B_2,
    B_3)$ and, at the same time, (1) are beyond either of the edges
    $(A_i, B_{i + 1})$ of $C_{A,B}$, and (2) are on the side of the
    plane $B_1 A_2 B_3$ opposite to $B_2$ or similarly for the
  analogous planes
    $B_1A_3B_2$, $B_2A_1B_3$ and the points  $B_3$, $B_1$
  respectively.

  \item The visible points beyond the facet $(A_1, A_2, A_3)$ of
    $C_{A,B}$ are exactly the points that are also in the visibility
    cone of the Sch\"onhardt polytope.
  \end{enumerate}
\end{lemma}

In what follows we will use the language of oriented matroids.
 For the theory of oriented matroids we refer to
\cite{RedBook} and \cite{Zie94}.  Here we only sketch the necessary
definitions and how they are related to the notion of visibility.  The
orientation of a simplex $(x_1, x_2, x_3, x_4)$, is defined as
$$
[x_1, x_2, x_3, x_4] = \mbox{sign } \det \left( 
  \begin{array}{cccc}
    x_1 & x_2 & x_3 & x_4 \\
    1 & 1 & 1 & 1
  \end{array}
\right).
$$
All such orientations make up the {\em chirotope} of an oriented matroid (see
page 123 in \cite{RedBook}).  

Given the oriented matroid of points $x_1, \ldots, x_n$ in $d$--space,
its circuits are functions $C : \{x_1, \ldots, x_n\} \mapsto \{+, -,
0\}$ that correspond to so--called minimal Radon partitions.  This
means that the convex hulls of $C^+ = \{x_i | C(x_i) = +\}$ and $C^- =
\{x_i | C(x_i) = -\}$ intersect in their relative interiors, and $C^+$
and $C^-$ are minimal at that.  It is easy to check that the function
$$
C(x) = \left\{ \begin{array}{ll} (-1)^i \cdot [\overbrace{x_1,
      \ldots, x_{d + 1}}^{\mbox{\small omit $x_i$}}] &
    \mbox{, if } x \in \{x_1, \ldots, x_{d + 1}\}, \\[.3cm]
    0 & \mbox{otherwise,}
  \end{array} \right.
$$
defines a circuit if it is not identical $0$. In fact, all circuits
can be obtained this way. We will compute
circuits to use a argument of the following form: $x$ does not see the
skylight if and only if there is a circuit such that the positive part
is one of the diagonals and negative part is the set containing $x$ and
a subset of vertices of the skylight. Since then the tetrahedron
spanned by $x$ and the skylight is pierced by the diagonal.

Important tools to compute simplex orientations are the
Grassmann--Pl\"ucker relations (see Section 2.4 in \cite{RedBook}):
For points $a$, $b$, $x_1, \ldots, x_4$ they state that the set of
signs

$$
\{[a, b, x_1, x_2]\cdot[a, b, x_3, x_4], -[a, b, x_1, x_3]\cdot[a,
b, x_2, x_4], [a, b, x_1, x_4]\cdot[a, b, x_2, x_3]\}
$$
is either identical $0$ or contains both a $+$ and a $-$.  The typical
use of the Grassmann--Pl\"ucker relations is to deduce one orientation
when the others are known. We can read the orientations of some of the
different tetrahedra from two-dimensional projections (drawings) of
the point configurations as in Figure \ref{Schonhardt}. We use a
left--handed rule system, i.e.~we decide whether the triangle $(x_1,
x_2, x_3)$ is oriented counterclockwise ($+$) or not ($-$), also if
$x_4$ is on our side of the plane spanned by $x_1$, $x_2$, and $x_3$
($+$) or not ($-$), and multiply these two signs to obtain the
orientation $[x_1, x_2, x_3, x_4]$.

\bigskip

\noindent {\em Proof of Lemma \ref{SchonhardtPolytopeLemma}:}

\noindent 1. In a Sch\"onhardt polytope, the simplices $(A_1, A_2,
A_3, B_1)$, and $(A_1, A_2, A_3, B_2)$ have the same orientation since
edges $(A_1, A_2)$ and $(A_2, A_3)$ are both incident to facet $(A_1,
A_2, A_3)$ and they are both non-reflex edges.

By the above argument, going around the boundary of a Sch\"onhardt
polytope, keeping in mind which edges are reflex, we can determine the
orientation of 12 simplices up to one global sign change (there are 12
edges).  But there are $\left(6 \atop 4\right) = 15$ simplices formed
by the vertices of the Sch\"onhardt polytope.  The remaining three
simplices are $(A_1, A_2, B_2, B_3)$, $(A_2,A_3,B_1,B_3)$,
$(A_1,A_3,B_1,B_2)$.  The signs are determined by the following
Grassmann-Pl\"ucker relations: For $(A_1, A_2, B_2, B_3)$ take $a =
A_1$, $b = A_2$, $x_1 = A_3$, $x_2 = B_1$, $x_3 = B_2$, $x_4 = B_3$
(the other two by circular index shift).  Then: 
$$
\{- \cdot ?, - \cdot - \cdot +, - \cdot -\} \supseteq \{+, -\}
$$
the equation forces $[A_1, A_2, B_2, B_3] = +$.
From the chirotope information it is easy to check that all vertices
are in convex position (see description of how to read the facets of
the convex hull from the chirotope in Chapter 3 of \cite{RedBook}), and
that their convex hull $C_{A,B}$ is indeed an octahedron.

\bigskip
\noindent 2. We will show that if a point $x$ lies beyond $A_1 B_2$
of $C_{A,B}$, on the side of $B_1 A_2 B_3$ opposite to $B_2$ , then
$(B_1, A_2)$ and the triangle $(B_2, B_3, x)$ form a minimal Radon
partition in the set of vertices $A_1$, $A_2$, $A_3$, $B_1$, $B_2$,
$B_3$, and $x$, hence have an interior point in common. This means $x$
cannot see the skylight. For this, we compute the following
orientations:

$$
\begin{array}{lcll}
  -[B_1, B_2, B_3, x] & = & +, & \mbox{since $(B_1, B_2, B_3)$ is a facet
    of $C_{A,B} \setminus (A_1, B_2)$,}  \\
  +[A_2, B_2, B_3, x] & = & +, & \mbox{since $(A_2, B_2, B_3)$ is a facet
    of $C_{A,B} \setminus (A_1, B_2)$,} \\
  -[A_2, B_1, B_3, x] & = & -, & \mbox{from the assumption on $x$,} \\
  +[A_2, B_1, B_2, x] & = & -, & \mbox{from the Grassmann--P\"ucker
    relation below,}\\
  -[A_2, B_1, B_2, B_3] & = & -, & \mbox{from Part (1).}
\end{array}
$$
The necessary Grassmann--Pl\"ucker relation is the one with $a =
B_1$, $b = B_2$, $x_1 = B_3$, $x_2 = A_1$, $x_3 = A_2$, and $x_4 = x$
such that
$$
\{- \cdot ?, - \cdot - \cdot -, - \cdot + \} \supseteq \{+, -\}
$$
forces $[B_1, B_2, A_2, x] = -$.
  

\bigskip
\noindent 3. If $x$ is in the visibility cone $V$, then it is, by part
(2) of this lemma, on the same side as $B_3$ with respect to the plane
$B_1 A_2 B_2$. Hence $A_2$ is on opposite side of $B_3$ with respect
to the plane $B_1 B_2 x$.  Therefore, the relative interior of the
convex hull of $B_1$ and $A_2$ lies strictly on one side of the plane
$B_1 B_2 x$, and the tetrahedron $(B_1, B_2, B_3, x)$ on the other
side of this plane.  Therefore those two point sets cannot have points
in common.  By symmetry it follows that the other two diagonals do not
obstruct any point of $V$ from seeing the skylight either.

Assume now that a point $x$ is beyond face $(A_1, A_2, A_3)$, but
outside of $V$, i.e.~for instance on the $A_1$ side of the plane $B_1
B_2 A_2$.  We claim that the pair $\{B_1, A_2\}, \{B_2, B_3, x\}$ forms
a circuit in the oriented matroid of the point configuration of the
vertices of $C_{A,B}$ and $x$.  This means that the triangle $(B_2, B_3, x)$
is pierced by the diagonal $(B_1, A_2)$ in the relative interior, hence
$x$ is not visible.
$$
\begin{array}{lcll}
  -[B_2, B_3, A_2, x] & = & -, & \mbox{since $(B_2, B_3, A_1)$ is a facet
    of $C_{A,B} \setminus (A_1, A_2, A_3)$,}  \\
  +[B_1, B_3, A_2, x] & = & +, & \mbox{from the Grassmann--Pl\"ucker
    relations below,} \\
  -[B_1, B_2, A_2, x] & = & +, & \mbox{from the assumption on $x$,} \\
  +[B_1, B_2, B_3, x] & = & -, & \mbox{since $(B_1, B_2, B_3)$ is a
    facet of $ C_{A,B}\setminus (A_1, A_2, A_3)$,}\\
  -[B_1, B_2, B_3, A_2] & = & +, & \mbox{from Part (1).}
\end{array}
$$ In this case, we have to apply the Grassmann--Pl\"ucker relations
twice to get $[B_1, B_3, A_2, x] = +$.  First we deduce $[A_1, A_2,
B_3, x] = -$ from the Grassmann--Pl\"ucker relation with $a = A_1$, $b
= A_2$, $x_1 = A_3$, $x_2 = x$, $x_3 = B_2$, $x_4 = B_3$: $$
\{ - \cdot -, - \cdot + \cdot ?, + \cdot +\} \supseteq \{+, -\}.
$$
Now we use this orientation to formulate $a = A_2$, $b = B_3$, $x_1
= A_1$, $x_2 = B_1$, $x_3 = B_2$, $x_4 = x$:
$$
\{ + \cdot -, - \cdot + \cdot ?, - \cdot +\} \supseteq \{+, -\}
$$
in order to get the desired $[A_2, B_3, B_1, x] = -$.
\endofproof

\subsection{The Cupola}

\begin{definition} \label{cupoladefinition} A polytope $C$ is called a
{\em cupola} if 
  it has the following properties:
  \begin{figure}[h]
    \psfrag{A1}{\small $A_1$}
    \psfrag{A2}{\small $A_2$}
    \psfrag{A3}{\small $A_3$}
    \psfrag{B1}{\small $B_1$}
    \psfrag{B2}{\small $B_2$}
    \psfrag{B3}{\small $B_3$}
    \psfrag{skylight}[l][l]{\small skylight}
    \psfrag{cupola C}[l][l]{\small cupola $C$}
    \psfrag{polytope P}[l][l]{\small polytope $P$}
    \begin{center} 
      \includegraphics[scale=.5]{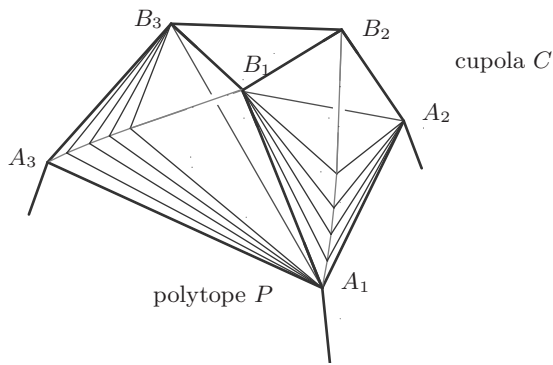}
    \end{center}
    \caption{A cupola as part of a larger convex polytope $P$}
    \label{cupola} 
  \end{figure}
  \begin{enumerate}
  \item the vertices of $C$ are $A_1, A_2, A_3$, $B_1, B_2, B_3$, and
    $q_k^{i, j}$ ($k = 0, \ldots, m+1$, $(i, j) \in \{(1, 2), (2, 3),
    (3, 1)\}$), where $q_0^{i, j} = A_i$ and $q_{m+1}^{i, j} = B_j$
    are identified.
  \item the vertices $A_1, A_2, A_3$, $B_1, B_2, B_3$ are in
    Sch\"onhardt position, and $(A_1, A_2, A_3)$ (the bottom facet)
    and $(B_1, B_2, B_3)$ (the skylight) are facets of $C$,

  \item the other facets are $(B_i, q_k^{i, j}, q_{k + 1}^{i, j})$ and
  $(A_j, q_k^{i, j}, q_{k + 1}^{i, j})$ for $k = 0, \ldots, m+1$, $(i,
  j) \in \{(1, 2), (2, 3), (3, 1)\}$, \item the vertices $q_k^{1, 2}$
  ($k = 1, \ldots, m$) lie on the side of the plane $B_1 A_2 B_3$
  opposite to $B_2$. Similar conditions must hold for $q_k^{2, 3}$ and
  $q_k^{3, 1}$.
  
  \end{enumerate}
\end{definition}

\begin{proposition} \label{cupolaprop} 
  Let $C$ be a cupola which is part of a larger polytope $P$, i.e.~$Q
  = P \setminus C$ is a convex polytope and $Q$ and $C$ share the
  common facet $(A_1, A_2, A_3)$.  Let $n$ be the number of vertices
  of $P$.
  
  If $T$ is a triangulation of $P$ with the property that the fourth
  point of the tetrahedron containing the skylight of $C$ is {\em not}
  in the visibility cone of $C$, then there are at least $n + m - 3$
  tetrahedra in the triangulation.
\end{proposition}

\proof If the vertex triangulating the skylight of $C$ is a vertex of
on a vertex-edge chain of $C$, then it does not see the skylight by
Definition \ref{cupoladefinition} (4) and Lemma
\ref{SchonhardtPolytopeLemma} (2).  If it is in $Q$ instead, then it
has to be beyond the face $(A_1, A_2, A_3)$ of $C$.  Hence by Lemma
\ref{SchonhardtPolytopeLemma} (3) it cannot see the skylight either.
Therefore the triangulation $T$ does not use one of the diagonals.  By
Lemma \ref{keyStructure} the number of tetrahedra is at least $n + m -
3$.  \endofproof


Later, Lemma \ref{openCondition} will be used to guarantee that we can
 place a point beyond a certain face.

\subsection{Constructing a Cupola from a Visibility Cone}

In this subsection we will show that cupolas can be attached to any
face of a frame polytope using intermediate polytopes and that the
visibility cone can be prescribed.  The following theorem does not
have the full strength we need for the construction.  In Section 3, we
will use a slightly stronger version which we will present at the end
of this section.  However, this theorem captures the main ideas used
to construct a cupola.

\begin{theorem} \label{cupolaConstruction} {\em (Cupola Construction
from a Given Visibility Cone)} Let $F$ be a facet of a $3$-polytope
$P$, and $V$ be a triangular cone such that $F \cap V$ is a triangle
in the relative interior of $F$, and $m$ be a positive integer.  Then
there is an $m$--cupola $C$ beyond $F$ of $P$ such that $P$ is beyond
$(A_1, A_2, A_3)$ of $C$ and such that $V$ is the visibility cone of
$C$.  Moreover, the input length of $C$ is polynomial in the input
lengths of $P$, $V$ and $m$.
\end{theorem}

Before we come to the proof, we will exhibit a necessary condition of
the visibility cone $V$ of a cupola $C$ and the facet the cupola is
being glued upon.  It will imply that we cannot directly attach a
cupola to a face (as in \cite{RupSei92}), but we have to construct an
intermediate polytope first.    
\begin{figure}[h]
  \psfrag{visibilityCone}[r][r]{\small visibility cone $V$}
  \begin{center}
    \includegraphics[scale=.5]{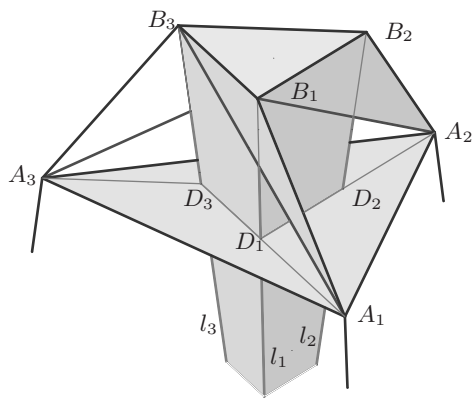}
  \end{center}
  \caption{Collinearity condition in the base triangle of a cupola}
  \label{necessaryCondition}
\end{figure}

\begin{lemma} \label{necessaryConditionLemma}  Let $A_1$, $A_2$,
  $A_3$, $B_1$, $B_2$, $B_3$ be vertices in Sch\"onhardt position.
  Define $l_1$ to be the intersection line of planes $B_3 B_1 A_1$ and
  $B_1 B_2 A_2$, lines $l_2$ and $l_3$ are defined accordingly
  (Figure~\ref{necessaryCondition}, note that they contain the extreme
  rays of $V$).  The lines $l_1$, $l_2$, and $l_3$ intersect the
  relative interior of the bottom face $(A_1, A_2, A_3)$ of a cupola
  $C$.  The intersection points $D_1$, $D_2$, and $D_3$ are forced to
  have the following collinearities: $A_1 D_1 D_2$, $A_2 D_2 D_3$, and
  $A_3 D_3 D_1$.
\end{lemma}

\proof $l_1$ enters the Sch\"onhardt polytope $S$ in point $B_1$, runs
along facet $(A_1, B_1, B_3)$ until it reaches the edge $(A_1, B_3)$
where it goes into the interior of $S$.  Then the relative interior of
$(A_1, A_2, A_3)$ contains the point $D_1$. In this way, $D_1$, $D_2$,
$A_2$ are all on the planes $A_1 A_2 A_3$ and $B_1 B_2
A_2$. \endofproof

\noindent {\em Proof of Theorem~\ref{cupolaConstruction}}: We proceed
in three steps.  The lines $l_1$, $l_2$, $l_3$ are defined as in Lemma
\ref{necessaryCondition}.  

\noindent {\em The bottom triangle $(A_1, A_2, A_3)$.}
We will now construct an intermediate polytope beyond $F$ which will
have a triangular facet $(A_1, A_2, A_3)$ which is (1) parallel to
$F$, and which is (2) intersected by the cone $V$ in a triangle $(D_1,
D_2, D_3)$ in the relative interior such that (3) the collinearity
condition from Lemma~\ref{necessaryConditionLemma} holds.  

\begin{figure}[h]
  \begin{center}
    \includegraphics[scale=.5]{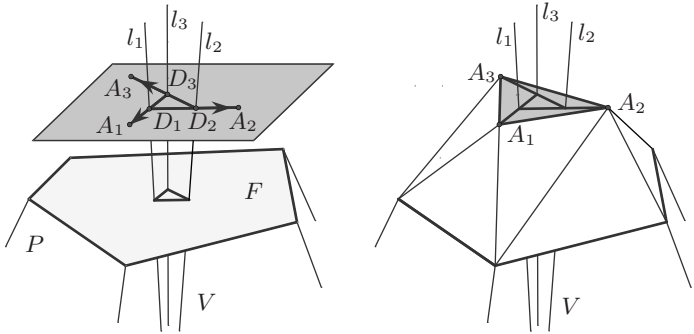}
  \end{center}
  \caption{Building the intermediate polytope for the cupola}
  \label{constructionSocket} 
\end{figure}

To do this, we place a plane $H$ parallel to and slightly above $F$
such that the intersection points $D_i$ of $H$ and $l_i$ ($i = 1, 2,
3$).  Also $H$ has to be so close to $F$ that the $l_i$ do not cross
between $H$ and $F$.  By prolonging the line segment $D_3 D_1$
slightly beyond $D_1$ (staying in $P \setminus F$) we obtain point
$A_1$, analogously construct $A_2$ and $A_3$ (Figure
\ref{constructionSocket}).  Taking the convex hull of $F$ and the
points $A_1$, $A_2$ and $A_3$ gives then the intermediate polytope,
whose face $(A_1, A_2, A_3)$ has the collinearity condition.  These
constructions are polynomially constructible in the sense of Lemma
\ref{openCondition}.

\begin{figure}[h]
  \psfrag{H'}{$H'$}
  \begin{center}
    \includegraphics[scale=.5]{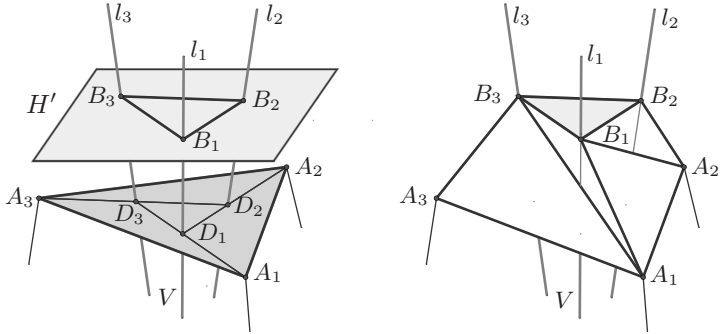}
  \end{center}
  \caption{Building the frame of a cupola} \label{buildframecup}
\end{figure}

\noindent {\em The frame of the cupola.} 
As in the construction of the bottom facet $(A_1, A_2, A_3)$, we place
a plane $H'$ parallel and slightly above this facet.  The intersection
of $H'$ and $V$ is the triangle $(B_1, B_2, B_3)$ ($B_1$ is on the
same extreme ray of $V$ as $D_1$ and so on). See figure
\ref{buildframecup}.
 
%

It is clear from the construction that triangles $(D_1, D_2, D_3)$,
$(A_1, A_2, A_3)$, and $(B_1, B_2, B_3)$ are parallel and all oriented
the same way.  Therefore it is not hard to check that the points
$A_1$, $A_2$, $A_3$, $B_1$, $B_2$, and $B_3$ are vertices of a
Sch\"onhardt polytope whose visibility cone is $V$.  Polynomiality of
this part of the construction follows from Lemma \ref{openCondition}
as well.

\noindent {\em Attaching the vertex-edge chains.}  Now that the frame
of a cupola done, i.e.~the vertices $A_1, \ldots, B_3$ in Sch\"onhardt
position, it remains to patch the key structures of
Lemma~\ref{keyStructure}, the vertex-edge chains $q_{i}^{j,k}$ ($i =
1, \ldots, m$, $(j,k) \in \{(1, 2), (2, 3), (3, 1)\}$), to the sides
of the frame $conv(P \cup \{A_1, \ldots B_3\})$.  

Given triangular faces $(a, q_0, q_{m+1})$ and $(b, q_0, q_{m+1})$ of
a convex polytope $P$ and a plane $G$ which (strictly) separates
points $q_0$ and $q_{m+1}$.  We claim that we can construct points
$q_1, \ldots, q_m$ beyond the edge $(q_0, q_{m+1})$ of $P$ such that
the convex hull of $P \cup \{q_i\}$ has the properties of
Lemma~\ref{keyStructure} and such that the points $q_1, \ldots, q_m$
lie on the same side of $G$ as $q_0$.  Moreover, the input length of
the constructed points is polynomially bounded in the input length of
$P$ and $G$.

By applying our claim three times, we will conclude our proof.  The
vertices $q_i^{j, j+1}$ are placed beyond edge $(A_j, B_{j+1})$,
vertices $B_j$ and $A_{j+1}$ take the roles of $a$ and $b$, $G$ is the
plane spanned by $B_j$, $A_{j+1}$ and $B_{j+2}$.  It is easy to check
that this is exactly what we want for Lemma~\ref{keyStructure} and for
the cupola conditions.

Now we prove the claim.  We will put the points $q_i$ ($i = 1, \ldots,
m$) on a parabola segment, beyond the edge $(q_0, q_{m+1})$.  Let $H$
be a plane containing $q_0$ and $q_{m+1}$ which also intersects the
interior of $P$.  This plane has the property that it contains points
beyond edge $(q_0, q_{m+1})$.  It is constructible in polynomial time.
(Let $v$ be the sum of the two normal vectors of planes $a q_0
q_{m+1}$ and $b q_0 q_{m+1}$, and $H$ the plane containing $q_0$ and
$q_{m+1}$ parallel to $v$.)

\begin{figure}[h]
  \psfrag{G}{\small $G$}
  \psfrag{a}{\small $a$}
  \psfrag{b}{\small $b$}
  \psfrag{q0}{\small $q_0$}
  \psfrag{q1}{\small $q_1$}
  \psfrag{qm}{\small $q_m$}
  \psfrag{qmp1}{\small $q_{m+1}$}
  \psfrag{E}{\small $E_\varepsilon$}
  \psfrag{ew}{\small $\varepsilon w$}
  \psfrag{D}{\small $D$}
  \begin{center}
    \includegraphics[scale=.6]{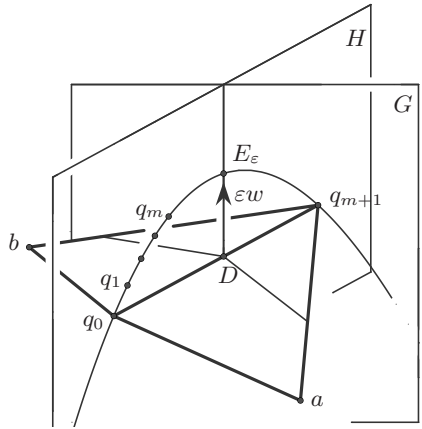}
  \end{center}
  \caption{Construction of the vertex-edge chain}
  \label{constructionVertexChain}
\end{figure}

Let now $D$ be the intersection point of $G$ and $(q_0, q_{m+1})$.
Let $w$ be a vector of direction of the intersection line of $G$ and
$H$, such that starting at $D$ it is pointing out of $P$.  Now let
$E_\varepsilon = D + \varepsilon w$ for $\varepsilon > 0$ to be
specified later.  For small $\varepsilon$, $E_\varepsilon$ is beyond
$(q_0, q_{m+1})$.  Hence the parabola defined according to Lemma
\ref{parabolaPlacement}, stated and proved below, by $p(0) = q_0$,
$p(1/2) = E_\epsilon$, and $p(1) = q_{m+1}$ lies entirely in $H$, and
for arguments between $0$ and $1$ passes just beyond $(q_0, q_{m+1})$.
Let $q_i = p(i / (4m))$ for $i = 1, \ldots, m$.  For small
$\varepsilon$ all those points are beyond $(q_0, q_{m+1})$ and on the
same side of $G$ as $q_0$ (polynomial conditions, use Lemma
\ref{openCondition}).  Also, they are in convex position such that the
convex hull of $P \cup \{q_1, \ldots, q_m\}$ has exactly the required
face lattice.  \endofproof

\begin{lemma} \label{parabolaPlacement} Let $p_0$, $p_1$,
  $p_2$ three non-collinear points in $\mathbb{R}^3$ and $t_0$, $t_1$,
  $t_2$ three distinct real numbers.  Then there is a unique curve $p:
  \mathbb{R} \to \mathbb{R}^3$ such that $p_0 = p(t_0)$, $p_1 =
  p(t_1)$, and $p_2 = p(t_2)$ which is quadratic in every coordinate.
  Furthermore, all points on $p(t)$ are in the plane spanned by $p_0$,
  $p_1$, and $p_2$, and they are in convex position.  Also a plane
  containing $p(r)$ and $p(l)$ for some $r \ne l$ which does not
  contain all of $p$ has all points between $l$ and $r$ on one of its
  sides and all other points on the other side.
\end{lemma}

\proof Since $p_0$, $p_1$, $p_2$ have to be on the $t_0$, $t_1$, $t_2$
positions of the curve
\[
p(t) = \left( \begin{array}{c}
    a_x + b_x t + c_x t^2 \\
    a_y + b_y t + c_y t^2 \\
    a_z + b_z t + c_z t^2 
  \end{array} \right), 
\]
we have the condition
\[
\left( \begin{array}{ccc} a_x & b_x & c_x \\ a_y & b_y & c_y \\ a_z &
    b_z & c_z \end{array} \right) 
\left( \begin{array}{ccc} 1 & 1 & 1 \\ t_0 & t_1 & t_2 \\ t_0^2 &
    t_1^2 & t_2^2 \end{array} \right) = 
\left( \begin{array}{ccc} \vdots & \vdots & \vdots \\ p_0 & p_1 & p_2
    \\ \vdots & \vdots & \vdots \end{array} \right).
\]
By the non-singularity of the Vandermonde matrices, there is a unique
solution to $a_.$, $b_.$, $c_.$ given the $p_i$ and $t_i$.
  
The curves which are quadratic in every coordinate are linear
transforms of the moment curve $m(t) = (1, t, t^2)$.  This curve lies
entirely in the $x = 1$ plane, is convex, and has the condition that
it intersected by each plane at most twice (or it is in this plane).
All these properties are invariant under linear transformations.
\endofproof

Proposition \ref{cupolaprop} stated that we get a large triangulation
if we triangulate the skylight of a cupola by a vertex outside the
visibility cone.  Now we want to estimate how much smaller a
triangulation is if we use a vertex $v$ in the visibility cone
instead.  We give a relatively small triangulation of the cupola and
of the space between the bottom face $(A_1, A_2, A_3)$ of the cupola
and the triangular face $F$ of $P$ with the help of the vertex $v$.

\begin{proposition} \label{triangulateCupola} Let $F$ be triangular
  face of a polytope $P$, and $C$ an $m$-cupola attached to it
  according to Lemma \ref{cupolaConstruction}.  Let $v$ be a vertex of
  $P$ in the visibility cone of $C$.  Then there is a triangulation of
  $\mbox{conv} (\{v\}, F, C)$ with at most $3m + 16$ tetrahedra.
\end{proposition}

\proof First of all, we triangulate along the vertex-edge chains using 
the tetrahedra $(B_i, A_{i + 1}, q_k^{i, i+1}, q_{k+1}^{i, i+ 1})$ for
$i = 1, 2, 3$, and $k = 0, \ldots, m$.

After removing these tetrahedra, we are left with the union of the
Sch\"onhardt polytope on the vertices $A_1$, $A_2$, $A_3$, $B_1$,
$B_2$, $B_3$, and the convex polytope $\mbox{conv} (\{v\}, F, (A_1,
A_2, A_3))$.  This is a non-convex polytope with all edges, except the
diagonals, being convex (easy conclusion from Lemma
\ref{SchonhardtPolytopeLemma} and the construction).  Since the
specified vertex $v$ is inside the visibility cone, it sees all facets
of this polytope, except the three facets it is incident to, from the
{\em interior}.  In particular, we can form tetrahedra of all these
facets and $v$ and none of them intersect badly.  They are $7$
tetrahedra for the facets of the Sch\"onhardt polytope (since we do
not count the bottom face) and at most 6 for the rest (the convex hull
of $F$ and $(A_1, A_2, A_3)$ has---by a planar graph argument---at
most $2 \cdot 6 - 4 = 8$ facets, subtracting 2 for $F$ and $(A_1,
A_2, A_3)$ gives 6).  \endofproof

It is this $3m$ in contrast to the $4m$ in Proposition
\ref{cupolaprop} which makes this triangulation better for large $m$.

\subsection{Constructing a Visibility Cone}

In order to use the cupola as a basic building block, we need to have
a visibility cone that contains a specified set of vertices and
intersects the relative interior of some face.  Once we have that we
can construct the cupola as described in the previous section.  The
set will consist of all vertices lying in a specified plane.

\begin{lemma} \label{oneWayToMakeVisibilityCones}
  Let $H$ be a plane which intersects the relative interior of some
  face $F$ of a polytope $P$, and let $S = \{v_1, \ldots, v_s\}$ be
  the set of vertices of $P$ lying in $H$, not including the vertices
  of $F$.  Let $S' = \{w_1, \ldots, w_{s'}\}$ a set of points in
  $\mbox{relint}(F) \cap H$.  It is possible to construct a triangular
  cone $V$ which intersects $F$ in a triangle that lies in the
  relative interior of $F$ and $V$ contains $S$ and $S'$ in its
  interior and no other vertex of $P$.
\end{lemma}

The reader may not see at this point the purpose of the set $S'$, but
we will justify it at the end of this section.

\begin{figure}[h]
  \psfrag{PcapH}{\small $P \cap H$}
  \psfrag{FcapH}{\small $F \cap H$}
  \psfrag{l}{\small $l$}
  \psfrag{r}{\small $r$}
  \psfrag{l'}{\small $l'$}
  \psfrag{r'}{\small $r'$}
  \psfrag{S'}{\small $S'$}
  \psfrag{S}{\small $S$}
  \psfrag{f'}{\small $f'$}
  \centerline{\includegraphics[scale=.5]{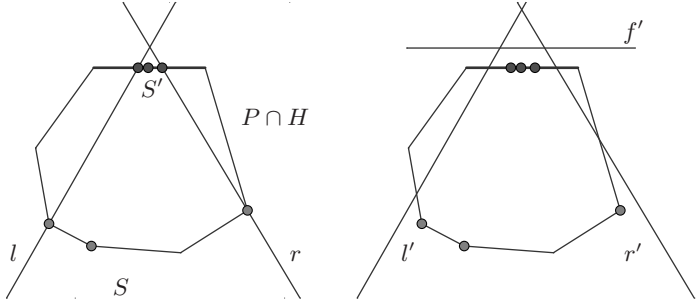}
    }
  \caption{Construction of $l$ and $r$, then $l'$ and $r'$ (viewed in $H$)}
  \label{visibilityConeConstruction1}
\end{figure}

\proof $P\cap H$ is a polygon.  Without loss of generality, $F \cap H$
is horizontal and situated on the top of the polygon $P \cap H$ (see
Figure \ref{visibilityConeConstruction1}).  Let $l$ be the line
connecting the leftmost point of $S'$ and leftmost vertex of $S$ (the
one encountered first when walking around $P \cap H$ counterclockwise,
starting at $F\cap H$).  Analogously, let $r$ be the line connecting
$M$ and the rightmost vertex of $S$.

The area between $l$ and $r$ (in $H$) is already a cone containing $S$
and no other vertices of $P$.  We will perturb it in a way that the
other conditions are satisfied as well.

First shift $l$ and $r$ parallely outwards, guaranteeing that they
still intersect $F\cap H$ in its relative interior (easy open
conditions); we obtain $l'$ and $r'$.  Also, let $f'$ be a line in $H$
parallel to $F$ just outside $P$, i.e.~such that $l'$ and $r'$
intersect $f'$ in the same order as $F \cap H$ (again using Lemma
\ref{openCondition}).

Now we will rotate $H$ about $l'$ and $r'$ and $f'$, getting three
planes bounding the desired triangular cone: Let $H$ be oriented in
some way, and $a_H x \ge b_H$ be its defining inequality.  Let $v$ be
some point which lies on the positive side of $H$.  Let $G_{l'}$ be
the plane through $l'$ and $v$.  By construction, all vertices in $S$
lie on the same side of $G_{l'}$, so we can orient it such that $S$ is
on its positive side.  Let $a_{l'} x \ge b_{l'}$ be its defining
inequality.  Perform the same construction for $r'$ and $f'$ obtaining
$G_{r'}$ and $G_{f'}$, also orienting them in a way that $v$ is on
their respective positive sides.  Let $G_{l'}^\varepsilon$ be the
plane defined by $(a_H + \varepsilon a_{l'}) x \ge b_{H} + \varepsilon
b_{l'}$.  This plane contains $l'$ and for small $\varepsilon$ it is
very close to $H$.  Hence, it is the {\em rotation of $H$ about $l'$
  in the direction of plane $G_{l'}$}.  Also let $G_{r'}^\varepsilon$
be defined by $(a_H + \varepsilon a_{r'}) x \ge b_H + \varepsilon
b_{r'}$, and $G_{f'}^\varepsilon$ be defined by $(-a_H + \varepsilon
a_{f'}) x \ge - b_H + \varepsilon b_{f'}$.

Obviously, all points in $S$ and in $S'$ are on the positive sides of
the planes $G_{l'}^\varepsilon$, $G_{r'}^\varepsilon$, and
$G_{f'}^\varepsilon$.  For small $\varepsilon > 0$, these planes do
not ``sweep'' over vertices of $P$ which are not in $S$, and it is
easy to see that in this case, there are no vertices of $P$ that
satisfy all three new inequalities.  Also for small $\varepsilon$, the
points in $F$ satisfying all three inequalities define a triangle in
the relative interior of $F$ with endpoints $G_{l'}^\varepsilon \cap
G_{r'}^\varepsilon \cap F$, $G_{r'}^\varepsilon \cap
G_{f'}^\varepsilon \cap F$, and $G_{f'}^\varepsilon \cap
G_{l'}^\varepsilon \cap F$.
\begin{figure}[h]
  \psfrag{Gle}{\small $G_{l'}^\epsilon$}
  \psfrag{Gre}{\small $G_{r'}^\epsilon$}
  \psfrag{Gfe}{\small $G_{f'}^\epsilon$}
  \psfrag{l'}{\small $l'$}
  \psfrag{r'}{\small $r'$}
  \psfrag{f'}{\small $f'$}
  \psfrag{v}{\small $v$}
  \psfrag{F}{\small $F$}
  \psfrag{Gl}{\small $G_{l'}$}
  \centerline{\includegraphics[scale=.6]{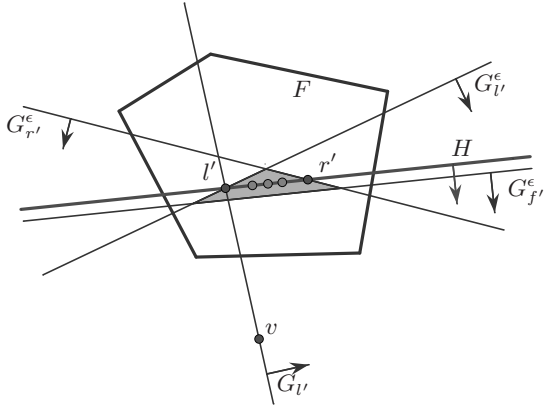}}
  \caption{Rotated hyperplanes, viewed by their intersections with $F$}
  \label{situationInF}
\end{figure}
Hence, the set of all points satisfying the three inequalities is a
triangular cone $V$ with the desired properties.  The conditions on
$\varepsilon$ are open polynomial conditions according to Lemma
\ref{openCondition}. \endofproof

This lemma can be used to build {\em one} cupola over the facet $F$.
However, there might be problems if we keep on constructing around the
polytope, like adding more cupolas over other facets of $P$.  The
visibility cone we just constructed might ``catch'' points we
construct later.  But these constructions all happen {\em beyond}
facets of $P$, so we can use the following lemma to construct all
cupolas one after the other without their visibility cones catching
extra vertices.

\begin{lemma}
  Let $H_1, \ldots, H_n$ hyperplanes, intersecting facets $F_1,
  \ldots, F_n$ of a polytope $P$ with the restriction that $F_i \cap
  H_j = \emptyset$ for all $i \neq j$.  Then $P_{\mbox{\scriptsize
      beyond }F_i} \cap H_j = \emptyset$ for all $i \neq j$.
\end{lemma}

\proof Assume there is a point $u$ in $P_{\mbox{\scriptsize beyond }
  F_i} \cap H_j$ ($i \neq j$).  Then this point also lies in
$(P\setminus F_i) \cap H_j$, but on the non-positive side of $F_i$.
Let $v$ be a point in $F_j \cap H_j$, Then $v$ is also in $(P\setminus
F_i) \cap H_j$ (since $F_j \subseteq P \subseteq P \setminus F_i$),
but on the positive side of $F_i$.  Hence, there must be a point $w$
on the line segment $[u,v]$ which is on the hyperplane containing
$F_i$.  The whole segment lies in $P \setminus F_i$, hence every point
on it has to satisfy all of $P$'s defining inequalities except that of
$F_i$.  So $w$ lies {\em in} the facet $F_i$.  But it also lies in
$H_j$ (the whole line segment does), which contradicts the assumption
$F_i \cap H_j = \emptyset$.
\endofproof

In Section 3.1 we will need an additional condition: Given a set of
lines in the plane $H$ (of Lemma \ref{oneWayToMakeVisibilityCones})
that pierce the face $F$, we want to be sure that these lines also
pierce the skylight of the constructed cupola.  (This condition will
play an important role when we want to force so-called blocking
conditions, see Section 3: At some point two tetrahedron spanned by
two skylights and two respective visible vertices $v$ and $v'$ will
have to intersect in their interiors.  This is already guaranteed if
the corresponding lines $g$ and $g'$ intersect inside the polytope.)

The next theorem specifies the way in which we will use all the
preceding lemmas in our construction in Section 3:

\begin{theorem} {\em (Full-Strength Cupola Construction)} \label{fullpower} 
  Let $H_i$ be planes that intersect facets $F_i$ of a polytope $P$ in
  their relative interiors such that $F_i \cap H_j = \emptyset$ for
  all $i \neq j$.  Let $S_i = \{v_1^i, \ldots, v_{s_i}^i\} := (vert(P)
  \cap H_i) \setminus F_i$, and $L_i = \{g_1^i, \ldots, g_{s_i'}^i\}$
  sets of lines.  Assume further that each of the lines $g_j^i$
  pierces the relative interior of $F_i$ and is incident to some
  $v_k^i$.
  
  Then we can sequentially construct all cupolas $C_i$ beyond the
  faces $F_i$ such that in the resulting polytope their visibility
  cones contain $S_i$ and no other vertices.  In addition, the
  skylight of the cupola over each $F_i$ is pierced by the lines in
  $L_i$.
\end{theorem} 

\proof The theorem follows from the ideas in Lemmas
\ref{cupolaConstruction} and \ref{oneWayToMakeVisibilityCones}.  In
the construction of the visibility cone over facet $F_i$, we invoke
Lemma \ref{oneWayToMakeVisibilityCones} with the polytope $P \cup
\bigcup_{j \ne i} P_{\mbox{\scriptsize beyond }F_j}$.  The set $S_i'$
is of course $\{l \cap F_i | l \in L_i\}$.  The cupola construction
was such that the cupolas over $F_j$ were always {\em beyond} the
facet $F_j$, so the constructed visibility cones contain no vertices
of the other cupolas.  In order to have the lines in $L_i$ pierce the
skylight of cupola $i$ we have to alter the construction of the cupola
in Lemma \ref{cupolaConstruction}: when we put the planes parallel to
$F$, we do it in such a way that the triangles $(A_1, A_2, A_3)$ and
then $(B_1, B_2, B_3)$ are pierced by these lines.  These are both
open conditions on the distance of the planes to $F$.\endofproof

\newcommand{\myfrac}[2]{{\displaystyle\frac{#1}{#2}}}

\section{The Transformation from SAT} \label{transformation}

It is our intention to model the well-known satisfiability problem
(SAT) using the visibility cones of cupola polytopes. Just as Ruppert
and Seidel did in \cite{RupSei92}, we will restrict our attention to
special SAT instances: each variable appears exactly three times,
twice unnegated and once negated.  This is not really necessary, but
simplifies explanations.  For our purpose this restriction is
appropriate because the SAT problem remains NP-complete even for
instances where each variable or its negation appear at most three
times (see references on page 259 in \cite{GarJohn79}).  In addition,
note that a change of variables can be used to change a non-negated
variable into a negated variable if necessary.  Also note that if a
variable appears only negated or only positive the variable and the
clauses that contain it can be discarded. Finally, if a variable
appears exactly once positive and exactly once negated then it can be
eliminated by combining the two clauses that contain the two variables
into one.  From now on, all logical formulas will have the properties
that each variable occurs exactly two times unnegated and exactly once
negated.  The formula
\[ 
f = (X_1 \lor \lnot X_2 \lor X_3 \lor \lnot X_4) \land (\lnot X_1 \lor
X_2 \lor \lnot X_3 \lor X_4) \land (X_1 \lor X_2 \lor X_3 \lor X_4)
\label{formula}
\]
is such a special SAT formula.  The figures in this section will
correspond to this particular instance.

In Section \ref{sectionLogical}, we will define the {\em logical
  polytope} associated to a given logical formula.  In Section
\ref{using} we will compute a number $K$ and see that the logical
polytope has a triangulation with $\le K$ tetrahedra if and only if
the logical formula is satisfiable.  Finally, in Section
\ref{sectionConstructing} we will give an algorithm to construct
explicit coordinates of a logical polytope.

\subsection{The Logical Polytope}
\label{sectionLogical}

\begin{figure}[h!]
  \center{
    \psfrag{variables}{4 variables} 
    \psfrag{clauses}{3 clauses}
    \psfrag{x}{\small$x$}
    \psfrag{y}{\small$y$}
    \psfrag{z}{\small$z$}
    \includegraphics[scale=.5]{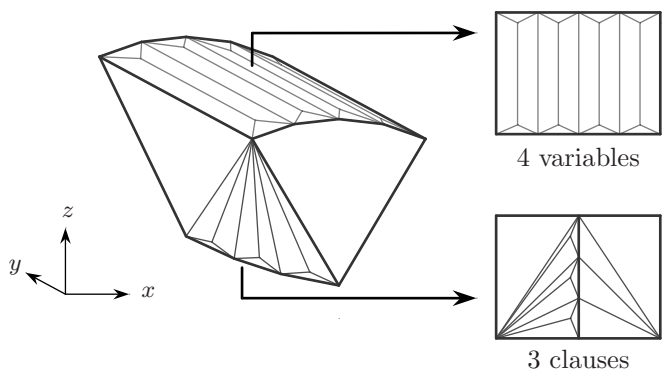}
    }
  \caption{Sketch of the logical polytope}
  \label{introframe}
\end{figure}

In the logical polytope there will be a cupola for each clause and one
for each variable and its negation. The cupolas will be glued to a
{\em frame polytope} which resembles a wedge.  Look carefully at
Figure \ref{introframe} for an example of the overall structure.

\begin{figure}[h!]
  \psfrag{clause1}{{clause 1}} \psfrag{clause2}{{clause 2}}
  \psfrag{clause3}{{clause 3}} \psfrag{c0}{\small$c_0$}
  \psfrag{c1}{\small$c_1$} \psfrag{c2}{\small$c_2$}
  \psfrag{c3}{\small$c_3$} \psfrag{c4}{\small$c_4$}
  \psfrag{c5}{\small$c_5$} \psfrag{c6=c2C}{\small$c_6 = c_{2C}$}
  \psfrag{x}{\small$x$} \psfrag{y}{\small$y$} \psfrag{z}{\small$z$}
  \centerline{ \includegraphics[scale=.5]{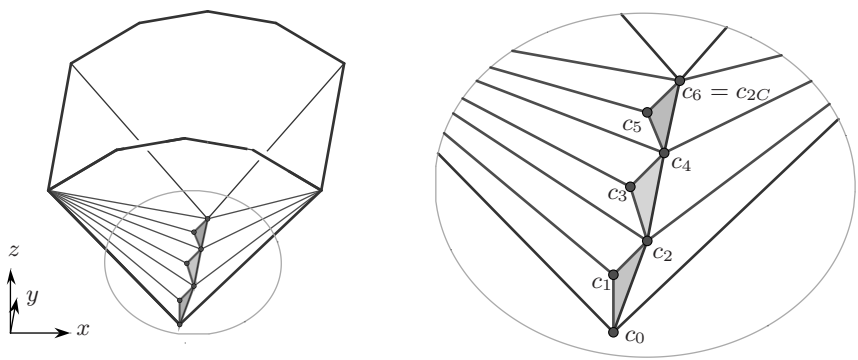} }
  \caption{The spine of the wedge: here the clause cupolas are
    attached}
  \label{belowwedge} 
\end{figure}

Figure \ref{belowwedge} displays the view of the lower hull of the
frame polytope, seen from the inside, i.e.~from above.  The sharp part
of the wedge consists of $2C+1$ vertices (where $C$ is the number of
clauses) $c_0,\dots,c_{2C}$.  We call this part of the frame polytope
the {\em spine}.  We attach the {\em clause cupola} associated with
clause $i$ to the triangle $(c_{2i},c_{2i+1},c_{2i+2})$ (shaded in the
picture).  

\begin{figure}[h!]
  \psfrag{zF1}{\small$z_F^1$}
  \psfrag{zF2=zT1}{\small$z_T^1 = z_F^2$}
  \psfrag{zR1}{\small$z_R^1$}
  \psfrag{zL1=zR2}{\small$z_L^1 = z_R^2$}
  \psfrag{zA1}{\small$z_A^1$}
  \psfrag{zC1}{\small$z_B^1$}
  \psfrag{zT4}{\small$z_T^V$}
  \psfrag{zL4}{\small$z_L^V$}
  \psfrag{variableCupola}{variable cupola}
  \centerline{
    \includegraphics[scale=.5]{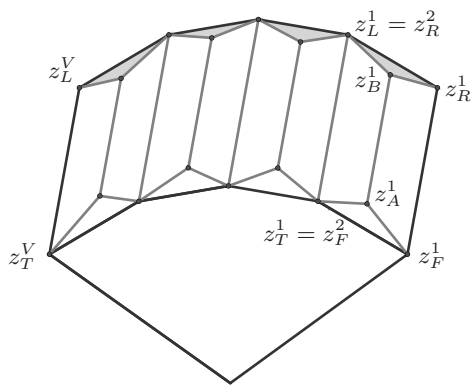}
    }
  \caption{The roofs, back gables shaded}  
  \label{roofs}
\end{figure}

On top of this wedge structure we will put a series of {\em roofs}.
They are triangular prisms, spanned by the two triangles $(x_T^i,
x_F^i, z_A^i)$ and $(x_L^i, z_R^i, z_B^i)$, one for every variable
$X_i$ of the logical formula. The {\em variable cupolas} will be
attached to the triangular facet $(z_L^i, z_R^i, z_B^i)$, the {\em
  back gables} (the triangular faces are shaded in Figure
\ref{roofs}).

The variable cupola of variable $X_i$ is such that its visibility cone
contains exactly the front vertices vertices $z_T^i$ and $z_F^i$. We
will use these cupolas to read from a small triangulation of the
polytope the logical value of variables according with the following
rule: if the truth-setting vertex $z_T^i$ associated to the $i$th
logical variable is used to triangulate the skylight of the cupola for
variable $i$, then we set $X_i = true$. If the truth-setting vertex
used to triangulate the skylight of the cupola for variable $i$ is
instead $z_F^i$ then $X_i = false$.

Beyond the quadrilateral face containing $z_T^i$ we will place the
{\em literal vertices} $x_1^i$ and $x_2^i$ which corresponds to the
positive occurrences of $X_i$ in the logical formula.  Beyond the
other quadrilateral face we will place the other literal vertex
$\overline{x_3^i}$ which correspond to the negated occurrence of
$X_i$.  These vertices are in the visibility cones of the three
cupolas of the clause where variable $X_i$ or its negation appears.

\begin{figure}[h!]
  \psfrag{x1}{\small$x_1^i$}
  \psfrag{x2}{\small$x_2^i$}
  \psfrag{x3}{\small$\overline{x_3^i}$}
  \psfrag{zF}{\small$z_F^i$}
  \psfrag{zT}{\small$z_T^i$}
  \psfrag{zR}{\small$z_R^i$}
  \psfrag{zL}{\small$z_L^i$}
  \psfrag{zA}{\small$z_A^i$}
  \psfrag{zC}{\small$z_B^i$}
  \psfrag{x}{\small$x$}
  \psfrag{y}{\small$y$}
  \psfrag{z}{\small$z$}
  \psfrag{variableCupola}{variable cupola}
  \centerline{
    \includegraphics[scale=.4]{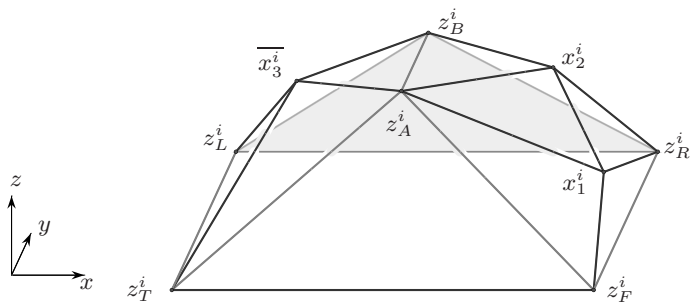}}
  \caption{A roof, back gable shaded, $z$-coordinate superelevated}  
  \label{roofcloseup}
\end{figure}

We list the five conditions on the logical polytope which are
necessary for the transformation to work in both ways, i.e.~a small
triangulation yields a satisfying truth assignment for our logical
formula and vice versa.

\begin{definition} For a given logical formula, a {\em logical
    polytope} is a three-dimensional polytope $P$ that satisfies the
  following conditions:
  \begin{enumerate} \label{logicalConstraints}
  \item {\em (Convexity)} The logical polytope must be convex and the
    face lattice is as we just described it.
  
\item {\em (Visibility)} The literal vertices $x^i_1$,
  $x_2^i$, and $\overline{x_3^i}$ are vertices in the visibility cone
  associated to their respective clause cupolas, but of no other
  clause visibility cone. The vertices $z_T^i$,$z_F^i$ are the only
  vertices in the visibility cones of the $i$th variable cupola.
  
\item {\em (Blocking)} This constraint ensures that the assignment of
  true or false values for variables is done consistently, i.e.~the
  positive (negative) literals can be used to make their clauses true
  if and only if the variable is set true (false).
  
  Concretely, the tetrahedron spanned by $z_F^i$ and the skylight of
  the cupola of variable $X_i$ intersects the interior of tetrahedron
  spanned by $x_1^i$ (by $x_2^i$) and the skylight of the clause
  cupola corresponding to $x_1^i$ (to $x_2^i$).  Also the tetrahedron
  spanned by $z_T^i$ and the skylight of the cupola of variable $X_i$
  intersects the interior of the tetrahedron spanned by
  $\overline{x_3^i}$ and the skylight of the clause cupola
  corresponding to it.  See Figure
\ref{blocking} for an example.  \begin{figure}[h!]  \centerline{
\psfrag{x1i}{\small$x_1^i$} \psfrag{zFi}{\small$z_F^i$}
\psfrag{zTi}{\small$z_T^i$} \psfrag{c_2j}{\small$c_{2j}$}
\psfrag{c2j+1}{\small$c_{2j+1}$} \psfrag{c2j+2}{\small$c_{2j+2}$}
\psfrag{skylight of variable Xi}{\small skylight of variable $X_i$}
\psfrag{skylight of clause j}{\small skylight of clause $j$}
\psfrag{roof of variable Xi}{\small roof of variable $X_i$}
\psfrag{x}{\small$x$} \psfrag{y}{\small$y$} \psfrag{z}{\small$z$}
\includegraphics[scale=.5]{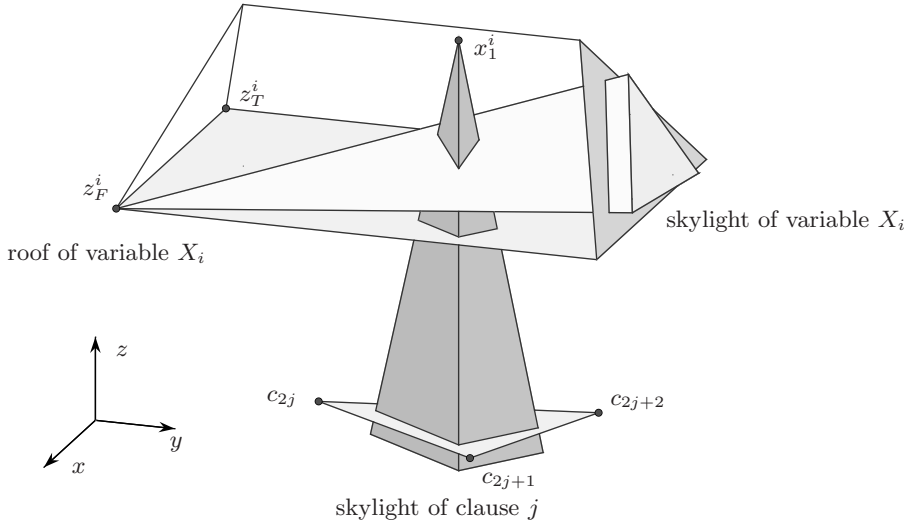}} \caption{Blocking for
consistent logical values} \label{blocking} \end{figure}
  
\item {\em (Non-blocking)} Using the vertex $z_T^i$ to triangulate the
  interior of the $i$-th variable cupola should not prevent the
  non-negated literal vertices from seeing their associated cupolas.
  Concretely, if $j$ is the clause corresponding to the literal vertex
  $x^i_1$, then tetrahedra $(z_T^i, z_L^i, z_R^i, z_B^i)$ and $(x^i_1,
  c_{2j-2}, c_{2j-1}, c_{2j})$ do not intersect at all.  The canonical
  analogue shall hold for $x^i_2$ and $\overline{x_3^i}$ (for
  $\bar{x_3^i}$ replace $z_T$ by $z_F$).
  
\item {\em (Sweeping)} We follow the same sweeping procedure proposed
  by Ruppert and Seidel \cite{RupSei92} we will need that
  
  (a) the variable $x^i_1$ is to the ``left'' (negative $x$ direction)
  of the planes $c_{2k-1}c_{2k}z_F^i$, $c_{2k}c_{2k+1}z_F^i$, and
  $c_{2k-1}c_{2k+1}z_F^i$ for $0 \leq k \leq C-1$.
  
  (b) $x_2^i$ is to the ``left'' of the planes $c_{2k-1}c_{2k}x^i_1$,
  $c_{2k}c_{2k+1}x^i_1$, and $c_{2k-1}c_{2k+1}x^i_1$ for $0 \leq k
  \leq C-1$.
  
  (c) $\overline{x_3^i}$ is to the ``left'' of the planes
  $c_{2k-1}c_{2k}z_F^i$, $c_{2k}c_{2k+1}z_F^i$, and
  $c_{2k-1}c_{2k+1}z_F^i$ for $0 \leq k \leq C-1$.
  
  (d) $z_T^i$ is to the ``left'' of the planes $c_{2k-1}c_{2k}x_2^i$,
  $c_{2k} c_{2k+1} x_2^i$, $c_{2k-1} c_{2k+1} x_2^i$, $c_{2k-1} c_{2k}
  \overline{x_3^i}$, $c_{2k} c_{2k+1} \overline{x_3^i}$, and $c_{2k-1}
  c_{2k+1} \overline{x_3^i}$ for $0 \leq k \leq C-1$.
\end{enumerate} 
\end{definition}

\subsection{Using the Logical Polytope}
\label{using}

\begin{lemma} \label{equivalence} Let $P$ be the Logical polytope, $m$
  the number of vertices on each vertex-edge chain, $n$ be the total
  number of its vertices.  For a SAT formula containing $C$ clauses on
  $V$ vertices there are polynomials $K(C, V)$ and $m(C, V)$ such that
  a logical polytope with $m = m(C,V)$ vertices on each vertex-edge
  chain admits a triangulation with $\le K = K(C, V) = n + m - 4$
  tetrahedra if and only if there is a satisfying truth assignment to
  the variables of the logical formula.
\end{lemma}

\proof If a triangulation $T$ of the polytope has $\le n + m - 4$
tetrahedra, then by Proposition \ref{cupolaprop} the skylight of each
cupola is triangulated by a vertex in the visibility cone of the
cupola.  In particular, one of $z_F^i$ and $z_T^i$ is chosen to
triangulate the cupola corresponding to variable $X_i$ for each $i$.
We claim that assigning to $X_i$ the truth value according to this
choice ($z_F \mapsto false$, $z_T \mapsto true$) satisfies all clauses
of the formula.

Each clause cupola skylight is triangulated by one of the literal
vertices, say clause $j$ by the positive literal vertex $x_1^i$ (or
$x_2^i$).  By the blocking conditions, it cannot be the case that the
variable skylight of $X_i$ is triangulated by $z_F^i$ since these
tetrahedra would intersect badly.  So we had set $X_i$ to $true$.
Having $x_1^i$ (or ${x_2^i}$) in clause $j$'s visibility cone meant
that variable $X_i$ appears unnegated in this clause.  If the skylight
of clause cupola $j$ is triangulated by $\overline{x_3^i}$, by the
same argument $X_i$ was set to false, and clause $j$ satisfied by the
literal $\lnot X_i$.  Hence all clauses are satisfied.

Now we need to prove the converse. If there is a $true$-$false$
assignment that satisfies all logical clauses we must find a
triangulation that has no more than $K$ tetrahedra.  For that we
construct a ``small'' triangulation.  There are three different kinds
of tetrahedra: the ones triangulating the cupolas, the ones
triangulating the roofs, and the ones triangulating of the rest of the
wedge.  We know how to triangulate a cupola if we know a vertex in its
visibility cone (see the proof of Prop.~\ref{triangulateCupola}).  For
the rest we will now follow a {\em sweeping procedure} which was first
described by Ruppert and Seidel \cite{RupSei92}.

The sweeping triangulation proceeds by triangulating ``slices'' that
correspond to the different variables $X_1$ to $X_V$, i.e.~from right
to left.  The variable roofs are arranged sequentially for exactly
this purpose.  A slice is roughly speaking the part of the tetrahedra
between a roof and vertices of the spine.  After the $i$th step of the
process the partial triangulation will have triangulated the first $i$
slices.  The part of the boundary of the partial triangulation that is
inside the logical polytope will form a triangulated disk.  We will
call it the {\em interface} following the convention of Ruppert and
Seidel.  It contains the following triangles:
\[
(z_T^i, c_{2C}, z_L^i) \mbox{ and } 
\left\{
  \begin{array}{p{3cm}cp{5cm}}
    $(z_T^i, c_{2j - 2}, c_{2j})$ & : & if clause $j$ is
      satisfied by one of the first $i$ variables, or \\
    $(z_T^i, c_{2j - 2}, c_{2j - 1})$ \\
    and  $(z_T^i, c_{2j -
      1}, c_{2j})$ & : & otherwise,
  \end{array}
\right.
\]
for all $j = 1, \ldots, C$.
\begin{figure}[h!]
  \psfrag{c0}{$c_0$}
  \psfrag{c2}{$c_2$}
  \psfrag{c3}{$c_3$}
  \psfrag{c4}{$c_4$}
  \psfrag{c5}{$c_5$}
  \psfrag{c2c}[r][r]{$c_6 = c_{2C}$}
  \psfrag{zT2}{$z_T^2 = z_F^3$}
  \psfrag{zL2}{$z_L^2 = z_R^3$}
  \centerline{\includegraphics[scale=.3]{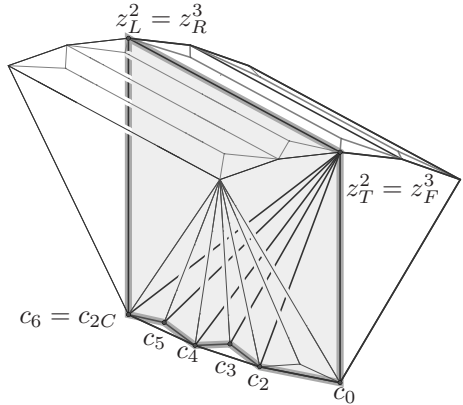}}
  \caption{The interface after step 2}
\end{figure}
Before the first step, the partial triangulation is empty.  After the
last step the partial triangulation will cover the whole logical
polytope.  In general, the vertices of the $i$th roof will see all
triangles of the interface and will be used as apexes to form new
tetrahedra to add to the current partial triangulation.  This way the
interface will slowly move from right to left. 

Now we describe in detail the triangulation step for the $i$th
variable $X_i$. Since we are only concerned with roof vertices in roof
$i$, we will drop all superscripts.  The triangulation step depends on
whether $X_i$ is set $true$ or $false$ in the satisfying assignment.
Let us consider first the case $X_i = true$:

\begin{figure}[h!]
  \psfrag{x1}{\small$x_1$}
  \psfrag{x2}{\small$x_2$}
  \psfrag{x3}{\small$\overline{x_3}$}
  \psfrag{zF}{\small$z_F$}
  \psfrag{zT}{\small$z_T$}
  \psfrag{zR}{\small$z_R$}
  \psfrag{zL}{\small$z_L$}
  \psfrag{zA}{\small$z_A$}
  \psfrag{zC}{\small$z_B$}
  \psfrag{x}{\small$x$}
  \psfrag{y}{\small$y$}
  \psfrag{z}{\small$z$}
  \psfrag{variableCupola}{variable cupola}
  \centerline{
    \includegraphics[scale=.45]{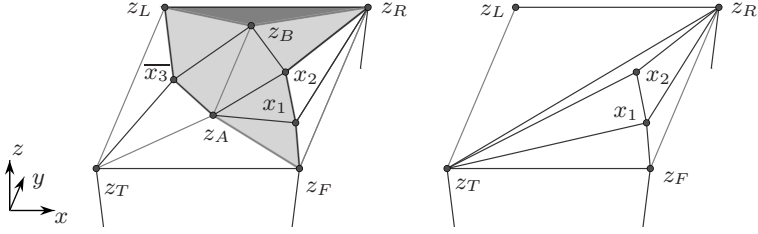}}
  \caption{Removing the tetrahedra spanned by $z_T$ and the shaded triangles}
  \label{removing}
\end{figure}

The point $z_T$ is used to triangulate the interior of the variable
cupola associated to $X_i$ according to Prop.~\ref{triangulateCupola}.
From $z_T$ we also form six tetrahedra with the following triangles:
$(z_L, \overline{x_3}, z_B)$, $(\overline{x_3}, z_B, z_A)$, $(z_B,
z_A, x_2)$, $(z_B, x_2, z_R)$, $(z_A, x_1, x_2)$, and $(x_1, z_A,
z_F)$. 

\begin{figure}[ht]
  \psfrag{x1}{\small$x_1$}
  \psfrag{x2}{\small$x_2$}
  \psfrag{x3}{\small$\overline{x_3}$}
  \psfrag{zF}{\small$z_F$}
  \psfrag{zT}{\small$z_T$}
  \psfrag{zR}{\small$z_R$}
  \psfrag{zL}{\small$z_L$}
  \psfrag{zA}{\small$z_A$}
  \psfrag{zC}{\small$z_B$}
  \psfrag{x}{\small$x$}
  \psfrag{y}{\small$y$}
  \psfrag{z}{\small$z$}
  \psfrag{c0}{\small$c_0$}
  \psfrag{c1}{\small$c_1$}
  \psfrag{c2}{\small$c_2$}
  \psfrag{spine}{\small spine}
  \psfrag{(a)}{ a.}
  \psfrag{(b)}{ b.}
  \psfrag{(c)}{ c.}
  \psfrag{(d)}{ d.}
  \psfrag{(e)}{ e.}
  \psfrag{(f)}{ f.}
  \centerline{
    \includegraphics[scale=.5]{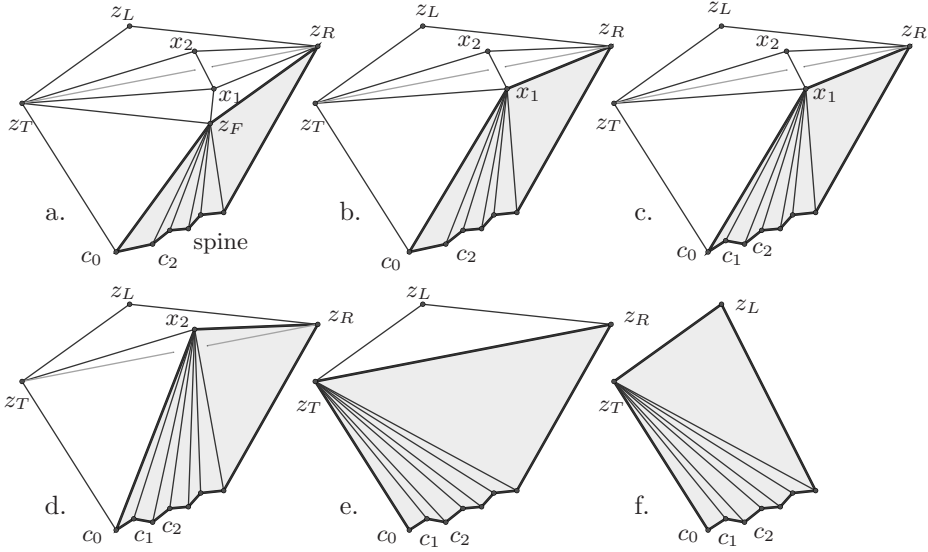}}
  \caption{The sweep}
  \label{sweep}
\end{figure}

Now we come to the part of the triangulation which gave the sweeping
procedure its name.  We form the tetrahedra between $x_1$ and the
current interface triangles. This is possible by part (a) of condition
5. We also use the tetrahedron $(x_1, z_T, c_{0}, z_F)$.  This is
illustrated in the transition from a.~to b.~in Figure \ref{sweep}.
The interface advances to $x_1$, i.e.~if $(z_F, c_j, c_k)$ was an
interface triangle before, now $(x_1, c_j, c_k)$ is an interface
triangle.  Also $(z_F, c_{2C}, z_R)$ is replaced by the triangle
$(x_1, c_{2C}, z_R)$.

Since $X_i$ is set to $true$ we can use $x_1$ to triangulate its
clause cupola according to Prop.~\ref{triangulateCupola}.  We only do
this if the clause cupola has not been previously triangulated using
an other literal vertex.  The condition 2 ensures that $x_1$ is in the
visibility cone of the clause cupola coming from the clause that
contains the unnegated literal $X_i$.  Furthermore, condition 4
ensures that we can actually perform this triangulation of the clause
cupola without badly intersecting the tetrahedra of the variable
cupola.  In Figure \ref{sweep}.c.~we see that if $x_1$ is used to
triangule clause $j$'s cupola, then the interface triangle $(x_1,
c_{2j-2}, c_{2j})$ is replaced by the two triangles $(x_1, c_{2j-2},
c_{2j-1})$ and $(x_1, c_{2j-1}, c_{2j})$.

We repeat this procedure with $x_2$, i.e.~form tetrahedra with $x_2$
and the current interface triangles, and then use $x_2$ to
triangulate its clause cupola if necessary (Figure \ref{sweep}.d.).
This is possible by part (b) of condition 5.  We continue by forming
tetrahedra using $z_T$ as apex (Figure \ref{sweep}.e, possible by
condition 5, part (d)).  At last, we will include the triangle
$(c_{2C}, z_L, z_B)$.  All these triangles are visible by part (d) of
the sixth condition.  After all these tetrahedra are added the
interface is ready for the next variable.

\begin{figure}[htb]
  \psfrag{x1}{\small$x_1$}
  \psfrag{x2}{\small$x_2$}
  \psfrag{x3}{\small$\overline{x_3}$}
  \psfrag{zF}{\small$z_F$}
  \psfrag{zT}{\small$z_T$}
  \psfrag{zR}{\small$z_R$}
  \psfrag{zL}{\small$z_L$}
  \psfrag{zA}{\small$z_A$}
  \psfrag{zC}{\small$z_B$}
  \psfrag{x}{\small$x$}
  \psfrag{y}{\small$y$}
  \psfrag{z}{\small$z$}
  \psfrag{c0}{\small$c_0$}
  \psfrag{c1}{\small$c_1$}
  \psfrag{c2}{\small$c_2$}
  \psfrag{spine}{\small spine}
  \psfrag{(a)}{ a.}
  \psfrag{(b)}{ b.}
  \psfrag{(c)}{ c.}
  \psfrag{(d)}{ d.}
  \psfrag{(e)}{ e.}
  \psfrag{(f)}{ f.}
  \centerline{
    \includegraphics[scale=.5]{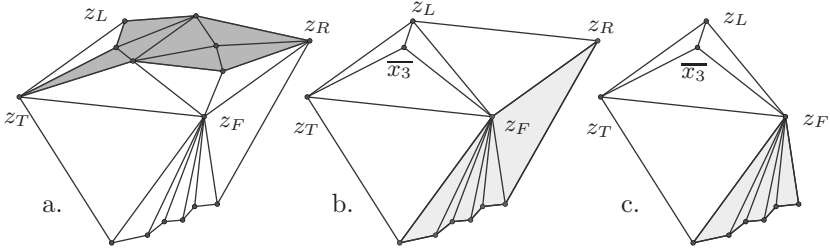}}
  \caption{The sweep for $X_i = false$}
  \label{sweepX3}
\end{figure}

Let us now consider the triangulation step in the case $X_i$ is set to
be $false$: We use the vertex $z_F$ to triangulate the variable cupola
as well as seven faces of the roof (see Figure \ref{sweepX3}): $(z_T,
\overline{x_3}, z_A)$, $(\overline{x_3}, z_A, z_B)$, $(\overline{x_3},
z_L, z_B)$, $(z_B, z_A, x_2)$, $(z_B, x_2, z_R)$, $(z_A, x_2, x_1)$,
$(x_2, x_1, z_R)$.  The reader can see that on the roof we are leaving
only the vertex $\overline{x_3}$.  Next the tetrahedron $(z_F, z_L,
z_R, c_{2C})$ is cut out.  Hereby the interface triangle $(z_F, z_R,
c_{2C})$ is replaced by $(z_F, z_L, c_{2C})$ (Figure
\ref{sweepX3}.c.).  Then $\overline{x_3}$ will be used as apex with
the triangles in the interface.  If the negated literal
$\overline{X_i}$ is used to satisfy its clause $j$, the $j$th clause
cupola is triangulated by $\overline{x_3}$.  The interface advances as
in the $true$-case.  Then $z_T$ can be used to form tetrahedra with
the triangles in the interface.  In the end the interface is again
ready for the next variable.

How may tetrahedra can such a triangulation have?  Triangulating all
cupolas with a vertex in their visibility cones yields at most $(3m +
16)(C + V)$ tetrahedra (Prop.~\ref{triangulateCupola}).  In one step of
the sweeping triangulation the tops of the roofs are each triangulated
using six or seven tetrahedra (if the variable is unnegated or
negated, resp.).  Furthermore, the interface is triangulated by some
vertices three times (in the positive case by $x_1^i$, by $x_2^i$, and
by $z_T^i$) or two times (in the negative case by $\overline{x_3^i}$
and by $z_T^i$).  The interface contains in each step between $C$ and
$2C$ triangles.  Eventually, in either case there is one more
tetrahedron (see above).  An upper bound for the size of this
triangulation is therefore
\begin{eqnarray*}
  \#T &\le& (3m + 16)(C + V) + 7V + 3 CV + 1 \\
  & = & m(3C + 3V) + \underbrace{16 C + 23 V + 3CV + 1}_{\displaystyle 
        p_T(V,C)} 
\end{eqnarray*}

What is the number of the vertices of the logical polytope in terms of
the number of clauses and variables?  We have $V$ logical variables
and $C$ clauses in the SAT instance. Say we have $m$ interior points
each of the vertex-edge chains we added (later we will replace the $m$
by a function of $V$ and $C$). We observe that we have $3m+6$ vertices
in each cupola, hence we have $(3m+6)(V+C)$ for all cupolas. We have
in each roof nine vertices, two of them are shared with the subsequent
roof except for the last roof.  Hence the total number of vertices in
roofs is $7V+2$. We have left only the $2C+1$ vertices along the
spine. In conclusion, the number of vertices of $P$ is
\begin{eqnarray*}
  n &=& (3m + 6)(V + C) + 7V + 2 + 2C + 1 \\
  &=& m (3C + 3V) + \underbrace{8C + 13V}_{\displaystyle p_n(C,V)} + 3 
\end{eqnarray*}

We had said before that a ``bad'' triangulation (where at least one
cupola skylight is triangulated by a vertex not lying in its
visibility cone) has at least $n + m - 3 = m (3C + 3V + 1) + p_n(C,
V)$ tetrahedra.  On the other hand a ``good'' triangulation has at
most $m (3C + 3V) + p_T(C, V)$ tetrahedra.  We can now set $m > p_n(C,
V) - p_T(C, V)$ and $K = m (3C + 3V) + p_T(C, V)$.  Then, if a good
triangulation exists, its size is smaller than or  equal to
$K$, and if not, all triangulations are larger than $K$.  Note finally
that the chosen $m$ and $K$ are polynomial in $C$ and $V$.
\endofproof

\begin{subsection}{Constructing the Logical Polytope}
\label{sectionConstructing}

\begin{lemma} \label{constructionPolytope} There is a polynomial
  algorithm that, given a logical formula on $V$ logical variables and
  $C$ logical clauses, produces a convex three-dimensional logical
  polytope as defined in Section \ref{sectionLogical}.  The
  coordinates of the vertices of the constructed polytope have binary
  encoding length polynomial in $V$ and $C$.
\end{lemma}

\proof The construction will be carried out in five stages. By the
time we end the construction all five requirements of the definition
of the logical polytope must be satisfied, but three of the conditions
will not be met until the last stage.

  \begin{enumerate}
  \item Give coordinates of the basic wedge, with rectangular faces
    on top for each variable.
  \item Attach the roofs for each variable, giving preliminary
    coordinates for the literal vertices and 
    preliminary coordinates for the points on the lower edge 
  (the spine of the wedge).
  \item Perturb the literal vertices to their final positions.
  \item Perturb the vertices on the spine of the wedge.
  \item Attaching the variable cupolas following the procedures of Section $2$.
  
  \end{enumerate}
  
  In every step we will build a construction element (a point, a
  line, or a plane) whose coordinates are polynomials in the
  construction elements up to that particular moment. Hence, the
  encoding length of each new construction element is bounded by a
  linear function of the encoding length of the construction so far.
  The number of construction steps is polynomially bounded in $C$ and
  $V$.  Hence the encoding length of the whole construction is also
  polynomially bounded in $C$ and $V$.  Note however, that the
  coordinates themselves will in general be exponentially large.
  
  Instead of writing explicit (and highly cumbersome) coordinates for
  the construction elements, we rely on Lemma \ref{openCondition} to
  ensure that such coordinates can be found if one has really the
  desire to see a particular logical polytope. A key property of
  Stages 2--4 in the construction is that the geometric conditions we
  want to determine a finite collection of strict polynomial
  inequalities {\em in a single variable}. Then, by Lemma
  \ref{openCondition}, we know there is an appropiate polynomial size
  solution. In subsequent stages of the construction similar new
  systems, for other independent parameters, will be solved,
  preserving what we had so far, but building up new properties.

\noindent  \emph{Stage 1: The basic wedge.}  
Consider the triangular prism which is the convex hull of the six
points $c_0 = (0, 0, 0)$, $c_{2C} = (0, 1, 0)$, $z_T^V = (0, 0, 1)$,
$z_F^1 = (1, 0, 1)$, $z_L^V = (0, 1, 1)$, and $z_R^1 = (1, 1, 1)$.
See Figure \ref{constructionWedge}(a).  In order to obtain a convex
structure on the top of the wedge, we consider the function $f(x) =
x(1-x) + 1$.  The vertices of each roof boundary (that is $z_T^i$ and
$z_F^i$ as well as $z_R^i$ and $z_L^i$) will lie on the surface $z =
f(x)$.  More specifically, $ z_F^i = z_T^{i+1} = (i / V , 0, f(i / V))
$ and $ z_R^i = y_L^{i + 1} = (i / V, 1, f(i / V))$ for $i = 0,
\ldots, n$.  By the concavity of $f$, the points are indeed in convex
position and their convex hull, the \emph{wedge} has the desired face
lattice (see Figure \ref{constructionWedge}(b)).

  \begin{figure}[ht]
    \psfrag{zLV}{\small $z_L^V$}
    \psfrag{zTV}{\small $z_T^V$}
    \psfrag{c0}{\small $c_0$}
    \psfrag{c2C}{\small $c_{2C}$}
    \psfrag{zF1}{\small $z_F^1$}
    \psfrag{zR1}{\small $z_R^1$}
    \psfrag{zR2}{\small $z_R^2 = z_L^1$}
    \psfrag{zF2}{\small $z_F^2 = z_T^1$}
    \psfrag{zA1}{\small $z_A^1$}
    \psfrag{zB1}{\small $z_B^1$}
    \psfrag{zA2}{\small $z_A^2$}
    \psfrag{zB2}{\small $z_B^2$}
    \psfrag{zA3}{\small $z_A^3$}
    \psfrag{zB3}{\small $z_B^3$}
    \psfrag{zA4}{\small $z_A^4$}
    \psfrag{zB4}{\small $z_B^4$}
    \psfrag{x}{\small $x$}
    \psfrag{y}{\small $y$}
    \psfrag{z}{\small $z$}
    \psfrag{a}{a.}
    \psfrag{b}{b.}
    \psfrag{c}{c.}
    \centerline{\includegraphics[scale=.5]{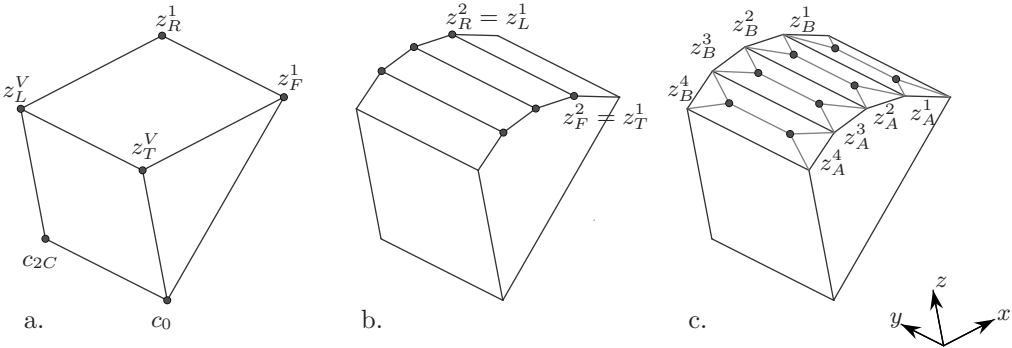}}
    \caption{Construction of the wedge}
    \label{constructionWedge}
  \end{figure}
  
  So far none of the conditions we want are satisfied (not even partially).

\noindent  \emph{Stage 2: The roofs.}  
We will first attach the points $z_A^i$ and $z_B^i$ to the
quadrilateral face $(z_L^i, z_R^i, z_T^i, z_F^i)$.  Then we give
preliminary coordinates to the literal vertices and to the vertices on
the spine.
  
Let $z_A^i = 1/2 \cdot (z_T^i + z_F^i) + (0, 1/3, t_{\mbox{\scriptsize
    roof}})$ and $z_B^i = 1/2 \cdot (z_T^i + z_F^i) + (0, 2/3,
t_{\mbox{\scriptsize roof}})$ where $t_{\mbox{\scriptsize roof}}$ is a
non-negative parameter that is called the {\em roof height}.  That is
the points have the same $x$ coordinate as the midpoint between
$z_T^i$ and $z_F^i$, $y$ coordinate $1/3$ and $2/3$ respectively, and
height $t_{\mbox{\scriptsize roof}}$ over the face $(z_T^i, z_F^i,
z_L^i, z_R^i)$.  We want to choose $t_{\mbox{\scriptsize roof}}$ in a
way that $z_A^i$ and $z_B^i$ are beyond the facet $(z_T^i, z_F^i,
z_L^i, z_R^i)$ (see Figure \ref{constructionWedge}(c)).  We can easily
achieve this by the technique presented in Lemma \ref{openCondition}:
The only possibly concave edges are the $(z_T^i, z_L^i)$.  One
restriction is therefore that all determinants $\det(z_T^i, z_L^i,
z_A^{i - 1}, z_A^i)$ have to be positive.  These are finitely many
open quadratic conditions on $t_{\mbox{\scriptsize roof}}$.  For
$t_{\mbox{\scriptsize roof}}=0$ the points $z_A^i$ an $z_B^i$ are
inside the facets $(z_T^i, z_F^i, z_L^i, z_R^i)$, hence the edges in
question are trivially convex.  We will get more polynomial
constraints on $t_{\mbox{\scriptsize roof}}$ below and then solve all
simultaneously to find the suitable roof height.

  The spine of the wedge is still a line.  We now put preliminary
  points $c_0, \ldots, c_{2C}$ on this line.  Let \[u(j) = \frac{1}{2}
  \frac{j}{2C}\] and $c_j = (0, u(j), 0)$ for $j = 0, \ldots, 2C - 1$,
  and $c_{2C} = (0, 1, 0)$ (see Figure \ref{spineLine}).  As an
  auxiliary point, let $b_l$ be the barycenter of the points $c_{2l -
  2}$, $c_{2l - 1}$, and $c_{2l}$ ($l = 1, \ldots, C$).  At this
  moment, this point $b_l=c_{2l -1}$.  Later, as we perturb the spine
  vertices $b_l$ will move accordingly, always $b_l = 1/3 (c_{2l - 2}
  + c_{2l - 1} + c_{2l})$.

  \begin{figure}[ht]
    \psfrag{c0}{\small $c_0$}
    \psfrag{c1}{\small $c_1$}
    \psfrag{c2}{\small $c_2$}
    \psfrag{c3}{\small $c_3$}
    \psfrag{c4}{\small $c_4$}
    \psfrag{c5}{\small $c_5$}
    \psfrag{c6}{\small $c_6$}
    \psfrag{c7}{\small $c_7$}
    \psfrag{c8}{\small $c_8$}
    \psfrag{0}{\small $0$}
    \psfrag{1/4}{\small $\displaystyle \frac{1}{4}$}
    \psfrag{1/2}{\small $\displaystyle \frac{1}{2}$}
    \psfrag{1}{\small $1$}
    \psfrag{x}{\small $x$}
    \psfrag{y}{\small $y$}
    \psfrag{z}{\small $z$}
    \centerline{\includegraphics[scale=.4]{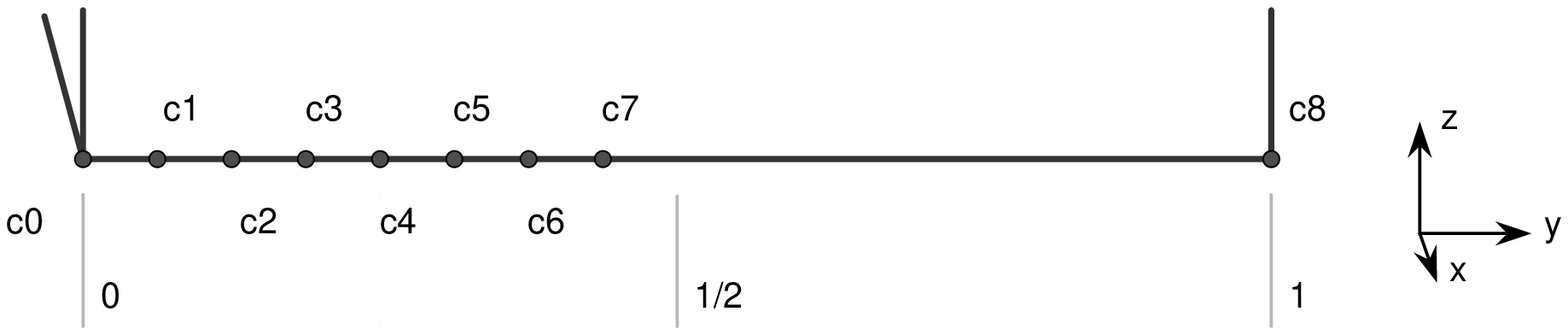}}
    \caption{Preliminary coordinates for the spine vertices}
    \label{spineLine}
  \end{figure}
  
  Now we want to give initial positions to the literal vertices.  Say
  variable $X_i$ occurs unnegated in clauses $l_1$ and $l_2$ and
  negated in $l_3$.  Note that $l_j$ depend on the variable we are
  considering.  For instance, in our example logical formula on
  p.~\pageref{formula}, for variable $X_1$, $l_1 = 1$, $l_2 = 3$, and
  $l_3 = 2$.  But for variable $X_2$, $l_1 = 2$, $l_2 = 3$, and $l_3 =
  1$.  
  
  The preliminary literal vertex $x_1^i$ is the intersection of the $y
  = u(2 l_1 - 1)$ plane with the line connecting $z_F^i$ and $z_B^i$.
  We do the same for the other positive occurrence of $X_i$ and obtain
  the preliminary $x_2^i$.  For the negative occurrence of $X_i$, we
  take the line connecting $z_T^i$ and $z_B^i$, intersect it with the
  $y = u(2 l_3 - 1)$ plane, and obtain the preliminary
  $\overline{x_3^i}$.  We join the preliminary $x_1^i$ and $b_{l_1}$
  by a line $d_{1}^i$ (this line lies in the $y = u(2 l_1 - 1)$
  plane).  Do the analogue process for $x_2^i$ and $\overline{x_3^i}$,
  obtaining $d_2^i$ and $d_3^i$.  Later we will move the vertices
  $x_1^i$, $x_2^i$, $\overline{x_3^i}$ along their respective lines
  $d_1^i$, $d_2^i$, $d_3^i$ a little out of polytope in order to turn
  them into extreme points. The lines $d_j^i$ will also be used for
  blocking conditions.

  \begin{figure}[hbt] 
    \psfrag{x1}{\small $x_1^i$} 
    \psfrag{x2}{\small $x_2^i$} 
    \psfrag{x3}{\small $\overline{x_3^i}$}
    \psfrag{cl1}[r][r]{\small $b_{l_1}$} 
    \psfrag{cl2}[r][r]{\small $b_{l_2}$} 
    \psfrag{cl3}[r][r]{\small $b_{l_3}$}
    \psfrag{zT}{\small $z_T^i$} 
    \psfrag{zF}{\small $z_F^i$}
    \psfrag{zC}{\small $z_B^i$} 
    \psfrag{zA}{\small $z_A^i$} 
    \psfrag{zR}{\small $z_R^i$} 
    \psfrag{zL}{\small $z_L^i$} 
    \psfrag{H}[r][r]{\small $y = u(2 l_1 - 1)$ plane} 
    \psfrag{a}{a.}  
    \psfrag{b}{b.}
    \psfrag{d1}{\small $d_1^i$}
    \centerline{\includegraphics[scale=.6]{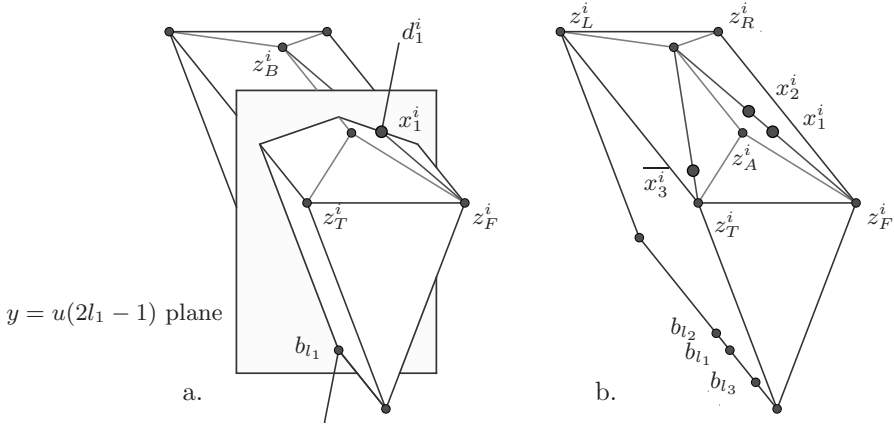}}
    \caption{Construction of the literal vertices in the $X_i$ slice of
       the wedge} 
    \label{constructionLiteral} 
  \end{figure} 
  
  Let $H^i$ be the plane that contains $z_T^i$ and $z_F^i$ and the
  midpoint of the edge $(z_L^i, z_B^i)$ (Figure \ref{constructionH1}).
  The only vertices above $H^i$ are $x_1^i$, $x_2^i$,
  $\overline{x_3^i}$, $z_A^i$, and $z_B^i$, and the only vertices {\em
  on} $H^i$ are $z_T^i$ and $z_F^i$.  This follows from the convexity
  of the current polytope.  

  \begin{figure}[ht]
    \psfrag{zF}{\small $z_F^i$}
    \psfrag{zT}{\small $z_T^i$}
    \psfrag{zA}{\small $z_A^i$}
    \psfrag{zC}{\small $z_B^i$}
    \psfrag{c0}{\small $c_0$}
    \psfrag{cl1}{\small $c_{2l_1 - 1}$}
    \psfrag{c2C}{\small $c_{2C}$}
    \psfrag{H}{\small $H^i$}
    \psfrag{g1}{\small $g_1^i$}
    \psfrag{x1}{\small $x_1^i$}
    \psfrag{d1}{\small $d_1^i$}
    \centerline{\includegraphics[scale=.5]{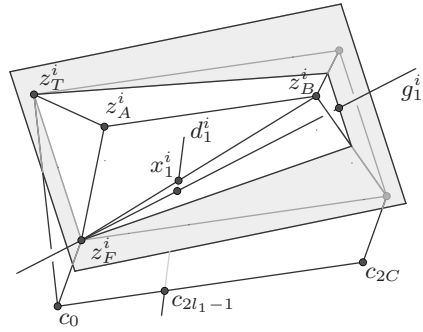}}
    \caption{Construction of $H^i$ and $g_1^i$}
    \label{constructionH1}
  \end{figure}
  
  Let $g_1^i$ ($g_2^i$) be the line in the plane $H^i$ which is
  incident to $z_F^i$ and intersects the line $d_1^i$ ($d_2^i$).  Note
  that this intersection point lies in the segment $(x_1^i, c_{2l_1 -
  1})$ (the line segment $(x_2^i, c_{2l_2 - 1})$), thus in the
  interior of the constructed polytope.  Analogously, let $g_3^i$ be
  the line in the plane $H^i$ which is incident to $z_T^i$ and
  intersects the line segment $(\overline{x_3^i}, c_{2l_3 - 1})$.  It
  can be verified that if the roof height is small $(z_L^i, z_R^i,
  z_B^i)$ is pierced by the $g_j^i$ in its relative interior. This is
  another strict polynomial inequality in $t_{\mbox{\scriptsize
  roof}}$.  It will be the planes $H^i$ and lines $g_j^i$ ($i = 1,
  \ldots, V$) from which we make the visibility cones for the cupolas
  of variables $X_i$ according to Theorem \ref{fullpower}.

  It is important to note right now that the non-blocking conditions
  are satisfied for this special position of the vertices.  We do not
  want the tetrahedron $(z_T^i, z_L^i, z_R^i, z_B^i)$ and the triangle
  $(x_1^i, c_{2l_1 - 2}, c_{2l_1})$ to intersect .  From this we get
  strict polynomial inequalities on $t_{\mbox{\scriptsize roof}}$.
  They are satisfied for $t_{\mbox{\scriptsize roof}} = 0$ since the
  $y$ coordinates of the spine vertices $c_l$ are smaller than $1/2$.
  A suitable value of $t_{\mbox{\scriptsize roof}}$ can be found
  solving the univariate inequality system we accumulated in our
  discussion (Lemma \ref{openCondition}).  It is easy to check that
  the sweeping conditions are also satisfied for the preliminary
  position of the points $x_1^i, x_2^i, \overline{x_3^i}$. So far we
  have met two of the five required conditions to have a logical
  polytope.

\noindent \emph{Stage 3: Literal vertices}
Now we put the final $x_j^i$ ($j = 1, 2, 3$) a little outward on line
$d_j^i$ (Figure~\ref{constructionLiteral}).  A little for $x_1^i$ and
$x_2^i$ means that the positive literal vertices lie in a plane
parallel to the face $(z_R^i, z_B^i, z_A^i, z_F^i)$ very close to it.
We treat $\overline{x_3^i}$ similarly. If the three literal vertices
are moved a sufficiently small distance $t_{\mbox{\scriptsize
    literal}}$ the face lattice of what we get after taking the convex
hull is as Figure \ref{roofcloseup} in all roofs.

By construction $H_i$ contains $z_F^i$ and $z_T^i$, and the $y = u(2j
- 1)$ planes contain all literal vertices corresponding to clause $j$.
This will become important for the visibility conditions (see Stage
5).  Also, for small $t_{\mbox{\scriptsize literal}}$ the non-blocking
and sweeping conditions are satisfied.
  
Although we do not have the blocking condition yet auxiliary lines can
be set up: As above, let $l_1$, $l_2$, $l_3$ be the clauses to which
the literal vertex $x_1^i$, $x_2^i$, $\overline{x_3^i}$ belong.  We
made sure that the line segments $(c_{2l_1-1}, x_1^i)$ and $(z_F^i,
z_B^i)$ intersect in their respective relative interiors.  Hence, by
the construction of line $g_1^i$, it is also pierced by $(x_1^i,
c_{2l_1 - 1})$ between $z_F$ and the face $(z_L^i, z_R^i, z_B^i)$.
(Analogously, $(c_{2l_2-1}, x_2^i)$ and $(z_F^i, g_2^i \cap (z_L^i,
z_R^i, z_B^i))$ as well as $(c_{2l_3-1}, \overline{x_3^i})$ and
$(z_T^i, g_3^i \cap (z_L^i, z_R^i, z_B^i))$ intersect in their
relative interiors). Later on this intersection will evolve into the
real blocking conditions using Theorem \ref{fullpower}.

\noindent \emph{Stage 4: The perturbing the vertices on the spine of the wedge.} 
We now perturb the points $c_j$ on the spine of the wedge.  Every
even-indexed $c_{2l}$ is changed to lie on a parabola, and for the
moment the odd-indexed vertices $c_{2l - 1}$ are changed to lie on the
line connecting $c_{2l - 2}$ and $c_{2l}$.  The $y$ coordinates of all
points stay the same:
\[
c_{2l} = \left(\frac{1}{2}(y - 1)^2 \cdot t_{\mbox{\scriptsize even}},
  y, (y - 1)^2 \cdot t_{\mbox{\scriptsize even}} \right).
\]  

Note that by the $1/2$ in the $x$ coordinate, the points are moved
\emph{into} the polytope.  The changes (parameter
$t_{\mbox{\scriptsize even}}$) must be small enough that the convex
hull now has the desired appearance (Figure \ref{perturbationSpine})
and the non-blocking conditions and the sweeping conditions are still
satisfied. Once more we appeal to Lemma \ref{openCondition}. The
polynomials inequalities are now on the variable $t_{\mbox{\scriptsize
    even}}$ and the sweeping and non-blocking were satisfied at
$t_{\mbox{\scriptsize even}} = 0$.  The reader should note that while
the constructed vertices in the roofs do not change coordinates,
dependent construction elements like the lines $d_j^i$ (connecting
$x_j^i$ and $c_{2l_j - 1}$) and $g_j^i$ (lying in $H^i$ and
intersecting $d_j^i$) change when the spine vertices move.  However,
the parameter $t_{\mbox{\scriptsize even}}$ has to be small enough
that the preliminary blocking conditions are still met: $g_j^i$ still
pierce the facet $(z_L^i, z_R^i, z_B^i)$ in its relative interior, and
$g_j^i$ and $d_j^i$ intersect in the interior of the polytope.
  
  \begin{figure}[ht]
    \psfrag{c0}{\small $c_0$}
    \psfrag{c1}{\small $c_1$}
    \psfrag{c2}{\small $c_2$}
    \psfrag{c3}{\small $c_3$}
    \psfrag{c4}{\small $c_4$}
    \psfrag{c5}{\small $c_5$}
    \psfrag{c6}{\small $c_6$}
    \centerline{\includegraphics[scale=.5]{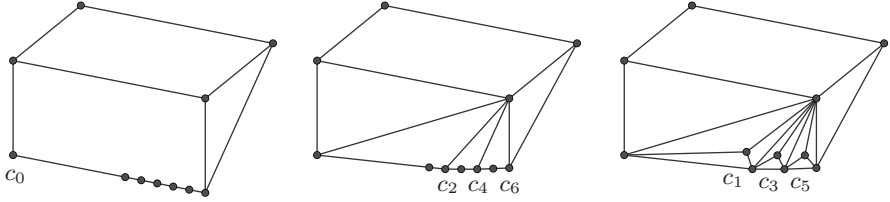}}
    \caption{Perturbation of the vertices on the spine}
    \label{perturbationSpine}
  \end{figure}
  
  Now we move the odd points $c_{2l - 1}$ beyond the face $G_l =
  (c_{2l - 2}, c_{2l}, z_T^0)$: to this end we choose a point $p_l$
  beyond $G_l$ and move to $c_{2l-1} + t_{\mbox{\scriptsize odd}} (p_l
  - c_{2l-1})$.  Such a point $p_l$ is easily found by taking a normal
  to $G_l$ through its barycenter and moving outwards while staying
  beyond the face (note that this involves again Lemma
  \ref{openCondition}, see the definition of {\em beyond}).  The
  parameter $t_{\mbox{\scriptsize odd}}$ is chosen small enough:
  Convexity and the correctness of the face lattice are easily
  achieved.  Also the sweeping conditions are valid for slight moves.
  Keeping $t_{\mbox{\scriptsize odd}}$ small also guarantees the
  non-blocking conditions: the tetrahedron $(x_1^i, c_{2l_1 -2},
  c_{2l_1 - 1}, c_{2l_1})$ is only slightly bigger than just the
  triangle $(x_1^i, c_{2l_1 -2}, c_{2l_1})$ which did not intersect
  the tetrahedron $(z_T^i, z_L^i, z_R^i, z_B^i)$ ($x_2^i$ and
  $\overline{x_3^i}$).
  
  For the blocking conditions, let $X_i$ be the $j$th logical variable
  in clause $l$.  Note that now the line $d_j^i$ intersects the
  triangle $(c_{2l - 2}, c_{2l - 1}, c_{2l})$ in its relative
  interior.  The lines $g_j^i$ are updated as the lines $d_j^i$ move.
  Since $t_{\mbox{\scriptsize odd}}$ is small, $g_j^i$ still pierces
  the facet $(z_L^i, z_R^i, z_B^i)$ in its relative interior, and
  $g_j^i$ and $d_j^i$ intersect in the interior of the polytope.  Note
  that $d_j^i$ is still in the $y = u(2l - 1)$ plane because the $y$
  coordinates of the spine vertices were conserved.
  
  \noindent \emph{Stage 5: Attaching the cupolas.}  It remains to
  construct all the cupolas.  Over the facets $(z_L^i, z_R^i, z_B^i)$
  ($i = 1, \ldots, V$) we construct cupolas using the planes $H_i$ and
  sets of lines $\{g_1^i, g_2^i, g_3^i\}$, and over the facets $(c_{2l
    - 2}, c_{2l - 1}, c_{2l})$ ($i = 1, \ldots, C$) we construct the
  clause cupolas using the $y = g(2l - 1)$ planes and the sets of
  lines $\{d_j^i \:|\: X_i\mbox{'s $j$th occurence is in clause }
  l\}$.  We invoke Theorem \ref{fullpower} and get the final polytope.
  By this construction, it is convex, has the correct face lattice,
  and the visibility conditions are satisfied.
  
  The reader will recall that $g_j^i$ and $d_j^i$ intersect in the
  interior of the polytope.  Say again variable $X_i$ occurs unnegated
  in clauses $l_1$ and $l_2$ and negated in $l_3$.  By Theorem
  \ref{fullpower} $g_j^i$ pierces the skylight of the cupola
  corresponding to variable $X_i$ and $d_j^i$ pierces the skylight
  corresponding to its clause $l_j$.  Hence, the tetrahedron spanned
  by $z_F^i$ and the variable $X_i$'s skylight and the tetrahedron
  spanned by $x_1^i$ ($x_2^i$) and clause $l_1$'s skylight ($l_2$'s
  skylight) intersect in their interiors.  Analogously, the
  tetrahedron spanned by $z_T^i$ and the variable $X_i$'s skylight and
  the tetrahedron spanned by $\overline{x_3^i}$ and clause $l_3$'s
  skylight intersect in their interiors.  These are exactly the
  blocking conditions.
  
  All other conditions concerned only points we constructed before, so
  they are still satisfied.  The final polytope is therefore a logical
  polytope.  \endofproof
  
  {\em Proof of Theorem \ref{main}:} The problem is clearly in $NP$:
  checking whether a collection of tetrahedra is indeed a
  triangulation of the polytope $P$ needs only a polynomial number of
  calculations.  Every pair of tetrahedra is checked for proper
  intersection (in a common face or not at all), and the sum of the
  volumes equals the volume of $P$ (computable for instance by the
  Delaunay triangulation of the polytope).  Also the size of
  triangulations of a given polytope is bounded by a polynomial in $n$
  of degree two (this follows from the well-known upper bound theorem,
  for details see \cite{rotstraus}).
  
  By Lemma \ref{constructionPolytope}, from a given logical formula on
  $V$ logical variables and $C$ clauses, we can construct a logical
  polytope $P$ of encoding length polynomial in $V$ and $C$.  Hence,
  by Lemma \ref{equivalence} there is a polynomial transformation that
  establishes the polynomial equivalence of a solution for the SAT
  problem and small triangulations of $P$. This completes the proof.
  \endofproof

\end{subsection}

\section{Related Problems and Conclusions} \label{conclusions}

In this last section we discuss Corollary \ref{maincoro} stated in the
introduction and relate our main theorem to previously known results
in the literature of optimal triangulations about the size of
triangulations. We conclude with some remarks.

{\em Proof of Corollary \ref{maincoro}:} For part (1), consider the
logical polytopes $P$ we have constructed.  If we could find the size
$s_{min}$ for their minimal triangulations in polynomial time, either
$s_{min} \leq K$, in which case we have indeed a triangulation smaller
or equal to $K$, or $ K < s_{min}$ in which case we can be sure there
is no triangulation for the logical polytope of size $K$ or less. This
proves it must be hard to find a minimal triangulation in dimension
$3$. To extend to any other dimension simply note that pyramids over
the logical polytopes have triangulations determined essentially by
the logical polytope.  In fact part (2) follows also from the pyramid
construction because the logical polytope is a face.  \endofproof

It is known that the sizes of triangulations for a $d$-dimensional
polytope with $n$ vertices lie between $n-d$ and $f_d(\partial
C(n+1,d+1))-d-1$, where $f_d(\partial C(n+1,d+1))$ is the number of
$d$-dimensional facets of a $d+1$ dimensional cyclic polytope with
$n+1$ vertices \cite{rotstraus}. In particular, for $3$-polytopes the
possible number of tetrahedra ranges from $n-3$ to ${n \choose
  2}-2n+3$.  Both bounds are known to be tight for three dimensions
\cite{edelsbrunner et al}. It is also known that the size of a minimal
triangulation of a convex $3$-polytope must lie between $n-3$ and
$2n-10$, when $n>12$ \cite{edelsbrunner et al}. That the upper bound
is tight was proved in \cite{sleatoretal} using hyperbolic geometry.
It is worth noticing at this point that the size of the constant $K$
we constructed in the previous section satisfies in fact $n-3 < K <
2n$.  Now we discuss an interesting reason why the lower bound is
strict:

We say that a convex polytope is {\em stacked} if it has a
triangulation whose dual graph is a tree (the dual graph of a
simplicial complex is the graph that has one vertex for each
maximal-dimensional simplex and two vertices are connected precisely
when the corresponding simplices are adjacent via a common facet). The
reader should be aware that in the literature the terminology stacked
polytope is often restricted to simplicial polytopes. Here of course,
we use it allowing that the stacking of simplices may give
coplanarities. For example, any $3$-cube or triangular prism is a
stacked polytope under our definition.

It turns out that a convex $d$-polytope $P$ with $n$ vertices has a
triangulation of size $n-d$ precisely when $P$ is a stacked polytope
(see \cite{rotstraus}).  It is natural to ask which polytopes are
stacked. If it were NP-hard to recognize stacked polytopes then this
would provide another proof that the problem of finding minimal
triangulations is also in the same class. However we can prove:

\begin{theorem}
  \begin{enumerate}
  \item For a convex $d$-polytope $P$ with $n$ vertices there is a
    $O(n^{d+2})$ algorithm to decided whether $P$ is stacked (i.e.,
    $P$ has a triangulation with $n-d$ maximal simplices). The
    algorithm uses only the 1-skeleton of the polytope. In particular,
    the size of a minimal triangulation of a stacked polytope does not
    depend on the particular coordinatization, but only depends on its
    face lattice.
  
  \item A convex $3$-dimensional polytope $P$ is stacked if and only
    if its graph does not contain as a minor the graph of an
    octahedron or a pentagonal prism.
  \end{enumerate}
\end{theorem}

\proof We need some definitions that have been introduced earlier in
the graph theory literature \cite{ArnborgProsku1,arncornpros,
  el-mallahcolbourn}. We say a graph $G$ is {\em $k$-decomposable} if
$G$ has $k+1$ or fewer vertices or there is a subset of vertices $S$
of $G$ with at most $k$ vertices such that 1) $S$ is a {\em cut}, i.e.
$G-S$ is disconnected, and 2) each of the connected components of
$G-S$ has the property that when the vertices of $S$ are added back
together with the complete graph on those vertices, the resulting
graph is again $k$-decomposable.

It was shown in Theorem 2.7 of \cite{ArnborgProsku1} that the class of
$k$-decomposable graphs is the same as the class of partial $k$-trees:
a graph is a {\em $k$-tree} if it can be reduced to the complete graph
$K_k$ , by a finite sequence of removals of degree $k$ vertices with
completely connected neighbors (i.e. neighbors of the vertex induce a
complete graph $K_k$). A {\em partial $k$-tree} is simply an
edge-subgraph of a $k$-tree.

Now we claim that a $d$-dimensional convex polytope $P$ is stacked, if
and only if its $1$-skeleton is a partial $d$-tree. Here is the proof:
the ``only if'' implication is clear from the definition of stacked
polytope.  We can prove the ``if'' implication by induction on $n$.
The theorem is trivial if $n=d+1$ because then $P$ is a simplex and
its graph is a complete graph, thus is a $d$-tree. Assume then $n>d+1$
and that the result is true for polytopes with fewer than $n$
vertices. Remember that if 1-skeleton $G(P)$ is a partial $d$-tree
then it is $d$-decomposable. Thus there is a cut $S$ of cardinality at
most $d$.  The set $S$ must have in fact cardinality $d$ because
$G(P)$ is $d$-connected by Balinski's theorem. If one has a
vertex-cutset $S$ of cardinality $d$ in $G(P)$, then the hyperplane
$H(S)$ spanned by $S$ intersects $G(P)$ only in the vertices of
$G(P)\cap S$ and in no edges (otherwise $S$ is not a cut).  In
conclusion, $H(S) \cap P$ is a $(d-1)$-simplex and because $G(P)$ is
$d$-decomposable we can apply induction hypothesis to prove the
polytopes $H(S)^+ \cap P$ and $H(S)^- \cap P$ are $d$-decomposable,
and thus they are partial $d$-trees with fewer vertices than $P$, so
both polytopes are stacked. Finally, note that their stacked
triangulations match well at the common boundary simplicial facet
$H(S) \cap P$, proving that $P$ is stacked.

Arnborg et al. discuss in \cite{arncornpros} an algorithm that, for
fixed values of $k$, decides whether a given graph is a partial
$k$-tree. The total running time is $O(n^{k+2})$. We also know that
computing the $1$-skeleton of a $k$-polytope can be done with the same
complexity. This, together with the above claim completes the proof of
the first part of the theorem.

For the second part we observe that partial $k$-trees form a minor
closed family. This means that the set of partial $k$-trees is closed
under taking edge-deletion or edge-contraction operations. The famous
results of \cite{robertsonseymour} imply that they can be
characterized by a finite set of forbidden minors. El-mallah and
Colbourn \cite{el-mallahcolbourn} proved that a graph is a planar
$3$-tree if and only if it has no minor isomorphic to the graph of an
octahedron or a pentagonal prism. This fact together with our claim
complete the proof of the second part.  \endofproof

\noindent {\bf Remarks:} 

1) The ``coning'' triangulation proposed in \cite{edelsbrunner et al}
provides an algorithm which is polynomial on the number of vertices
and gives a $2$-approximation of the minimal triangulation as it
produces a triangulation of size less than or equal to $2n-7$.

2) Given a $3$-dimensional convex polytope $P$, a proper subset $S$ of
tetrahedra with vertices in $vertices(P)$, and a positive integer
$K$. Deciding whether there is a triangulation of $P$ that uses $K$
simplices from $S$ can also be proved  to be an NP-hard problem using 
the constructions we explained.

3) It is interesting to note that the constructions presented in
\cite{Below99} prove also that {\em covering} a convex $3$-polytope
with tetrahedra can be done with fewer pieces than triangulating. A
cover is a collection of simplices whose union is the whole polytope,
but the elements can intersect in their interiors. It is would be
interesting to know what is the computational complexity of finding
minimal simplicial covers. Another interesting question would be what
is the complexity of deciding whether a triangulation of the boundary
of convex non-simplicial $3$-polytope extends to a triangulation of
the whole polytope without adding new interior points?  If hard, this
could be used to establish another proof our results.  The curious
reader can easily prove that already for a triangular prism not all
triangulations of the boundary extend to a triangulation of the whole
polytope.  Even more interesting. The triangulations of the boundary
of a $3$-cube extend or not depending on the coordinates of its
vertices.

\end{document}